\documentclass[11pt]{amsart}
\usepackage{geometry}  
\usepackage{amsthm}
\usepackage{url}
\usepackage{graphicx}
\usepackage{amssymb}
\usepackage{epstopdf}
\usepackage{color}
\usepackage{amsmath}
\usepackage{upgreek}
\usepackage{url}
\usepackage[square,numbers]{natbib}
\usepackage{subeqnarray}
\definecolor{darkred}{rgb}{.6,0,0}
\definecolor{darkblue}{rgb}{0,0,.7}

\newtheorem{lemma}{Lemma}[section]
\newtheorem{example}{Example}[section]
\newtheorem{remark}{Remark}[section]

\newtheorem{theorem}{Theorem}[section]
\newtheorem{assumption}{Assumption}[section]

\title[On a mean-field Pontryagin minimum principle]{On a mean-field Pontryagin minimum principle for stochastic optimal control}
\author{Manfred Opper}
\address{Institut f\"ur Softwaretechnik und Theoretische Informatik, Technische Universit\"at Berlin, Marchstraße 23, 10587 Berlin}
\author{Sebastian Reich}
\address{Institut f\"ur Mathematik, Universit\"at Potsdam, Karl-Liebknecht-Str. 24/25, 14476 Potsdam}
\date{\today}                                           

\begin{document}

\begin{abstract} 
This paper outlines a novel extension of the classical Pontryagin minimum (maximum) principle to stochastic optimal control problems. Contrary to the well-known stochastic Pontryagin minimum principle involving forward-backward stochastic differential equations, the proposed formulation is deterministic and of mean-field type. We denote it by the McKean--Pontryagin minimum principle. The Hamiltonian structure of the proposed McKean--Pontryagin minimum principle is achieved via the introduction of a pair of auxiliary functions. A gauge freedom in the choice of one of these two functions can be used to decouple the forward and reverse time equations; hence simplifying the solution of the underlying boundary value problem. We also consider  infinite horizon discounted cost optimal control problems. In this case, the mean-field formulation allows one to convert the computation of the desired optimal control law into solving a pair of forward mean-field ordinary differential equations. The McKean--Pontryagin minimum principle is tested numerically for a controlled diffusion process in a double well, a controlled inverted pendulum, a controlled Lorenz-63 system, and a controlled Lorenz-96 system. Although the focus is on linear-quadratic control problems, the proposed methodology is extendable to more general problems including mean-field type control formulations. 
\end{abstract}

\maketitle


%
\section{Introduction}
%

We consider the optimal control problem for a controlled stochastic differential equation (SDE)
of the form
\begin{equation} \label{eq:FSDE}
{\rm d}{\rm X}_t = b({\rm X}_t){\rm d}t + G({\rm X}_t) U_t {\rm d}t + \Sigma^{1/2} {\rm d}B_t, \qquad {\rm X}_0 = x_0,
\end{equation}
under either the finite horizon cost function
\begin{equation} \label{eq:costJ}
J_T(x_0,U_{0:T}) = \mathbb{E} \left[ \int_0^T \left(c({\rm X}_t) + \frac{1}{2} U_t^{\rm T} R^{-1} U_t \right) {\rm d}t + f({\rm X}_T) \right]
\end{equation}
or the infinite horizon discounted cost function
\begin{equation} \label{eq:costJ infinite}
    J(x_0,U_{0:\infty}) = \mathbb{E}\left[\int_0^\infty e^{-\gamma t}\left(c({\rm X}_t)+ \frac{1}{2} U_t^{\rm T} R^{-1} U_t\right){\rm d}t \right].
\end{equation}
Here $B_t$ denotes $d_x$-dimensional Brownian motion, $\Sigma \in \mathbb{R}^{d_x\times d_x}$ the symmetric positive semi-definite diffusion matrix, $R \in \mathbb{R}^{d_u\times d_u}$ a symmetric positive definite weight matrix, $G(x) \in \mathbb{R}^{d_x\times d_u}$ the possibly position-dependent control matrix, $c(x)\ge 0$ the running cost, $f(x)\ge0$ the terminal cost, $x_0 \in \mathbb{R}^{d_x}$ the initial condition, and $\gamma >0$ the discount factor. See, for example, \cite{Carmona,Meyn} for more details. We also introduce the weighted norm $\| \cdot \|_R$ via
\begin{equation} \label{eq:norm}
\|u\|^2_R = u^{\rm T} R^{-1} u.
\end{equation}
The aim is to find the closed loop control law $U_t = u_t({\rm X}_t)$ that minimizes the cost (\ref{eq:costJ}) or (\ref{eq:costJ infinite}), respectively, over the set of admissible control laws.

It is well known that the associated optimal control problems for ordinary differential equations (ODEs); {\it i.e.}~$\Sigma = 0$ in (\ref{eq:FSDE}), can be solved by the Pontryagin minimum (maximum) principle (PMP) \cite{pontryagin}. The PMP has been extended to controlled SDEs (\ref{eq:FSDE}) in terms of forward-backward SDEs giving rise to the so-called stochastic PMP \cite{Peng1993,Carmona,Bensoussan}. In this paper, we provide an alternative extension via a deterministic mean-field formulation, which leads to interacting particle implementations in the form of Hamiltonian ODEs. The formulation is attractive from a computational perspective since it allows for a decoupling of the forward and backward evolution equations for finite horizon optimal control problems and a forward simulation for infinite horizon optimal control problems. Our mean-field approach delivers closed loop control laws in contrast to the classical PMP for ODEs \cite{pontryagin}. The proposed methodology can then also be used to approximate viscosity solutions for controlled ODEs by setting $\Sigma = \epsilon I$ for sufficiently small $\epsilon >0$ \cite{CL83}.

Interacting particle approximations and their mean-field limits have been considered before in the context of the Hamilton--Jacobi--Bellman (HJB) equation for the value function $v_t(x)$ \cite{Carmona,Bensoussan}. The key idea is to connect the time evolution of the value function to a probability density function evolving under an appropriate Fokker--Planck equation and its mean-field representation \cite{MO22,JTMM22,reich23,Prashant24}. Alternatively, \cite{ZOL24} considers a coupled system of mean-field ODEs in the state variable $X_t$ and the value function $Y_t = v_t(X_t)$. In this paper, we directly approximate the gradient of the value function and derive appropriate Hamiltonian evolution equations; thus extending mean-field formulations to the PMP. We denote our mean-field formulation by the McKean--Pontryagin minimum principle. 

As a further alternative, we mention the application of direct optimization methods to learn a parametric form of the desired control laws in the form of $u_t(x;\theta)$ for finite horizon cost functions (\ref{eq:costJ}). See, for example, \cite{EHJ17,EHJ21,nusken2020solving,hua2025an,LVR24}, which rely on deep learning approximation architectures. The methods differ in the underlying procedure for sampling along a stochastic  forward process and in the definition of the cost function to be minimized. Similarly, infinite horizon discounted cost problems  can also be reformulated as  optimization problems in parameters $\theta$ based on approximations of the underlying stationary HJB equation; an approach pioneered in reinforcement learning \cite{Meyn}. 

We finally mention the saddle point formulation of \cite{MLLO24}, which is closely related to the Eulerian variational formulation stated in Section \ref{sec:Hamiltonian} and which can be solved by available algorithms for saddle point problems. See Remark \ref{rem1} for more details. Our approach is to instead transfer to a corresponding Lagrangian variational formulation \cite{Salmon88}, which then leads to our McKean--Pontryagin minimum principle.

Although the focus of this paper is on standard stochastic optimal control problems, our PMP formulation can be extended to mean-field control problems as described and analyzed, for example, in \cite{LI2012366,BLM16,Carmona}. Mean-field control problems arise in particular from partially observed diffusion processes \cite{Bensoussan92}. See \cite{R25c} for a first application of the McKean--Pontryagin minimum principle in combination with the ensemble Kalman filter for data assimilation \cite{reich2015probabilistic,law2015data,Evensenetal2022}, where mean-field variables appear in both the cost function (\ref{eq:costJ}) and the forward SDE (\ref{eq:FSDE}).

The remainder of the paper is structured as follows. The following Section \ref{sec:math} summarizes key mathematical results for finite horizon optimal control problems at the level of the underlying HJB and Fokker--Planck equations \cite{Carmona,Pavliotis2016}. The novel mean-field formulation of the PMP is provided in Section \ref{sec:Hamiltonian} starting from an Eulerian perspective in terms of suitable Hamiltonian partial differential equations (PDEs). Following \cite{Salmon88,holmI}, the arising Hamiltonian functional is subsequently transformed into a Lagrangian formulation with trajectories following mean-field Hamiltonian ODEs. The infinite horizon optimal control problem is treated in Section \ref{sec:discounted}. Here the key result is a pair of forward mean-field ODEs, which asymptotically deliver the desired closed loop control law. The proposed formulations are discussed in the context of linear control problems in Section \ref{sec:linear}. Several numerical approximation techniques are discussed in Section \ref{sec:implementation}; including an approximation based on discrete Schr\"odinger bridges \cite{MarshallCoifman,WR20} as outlined in Section \ref{sec:DM}. The McKean--Pontryagin formulation for finite horizon cost
(\ref{eq:costJ}) is carefully tested and compared to a numerical approximation of the HJB equation for a one-dimensional controlled diffusion process in Section \ref{sec:experiments}. Further numerical experiments demonstrate the broad scope of the proposed methodology, in particular for infinite horizon discounted cost (\ref{eq:costJ infinite}). The first example is provided by a controlled inverted pendulum \cite{Meyn} while the second example concerns a controlled Lorenz-63 system \cite{KK24}. The third numerical example in Section \ref{sec:num_IP2} presents results from an implementation of the McKean--Pontryagin formulation in a receding horizon setting applied again to a controlled inverted pendulum. The final example demonstrates the applicability of our approach to higher dimensional control problems by means of a controlled Lorenz-96 system \cite{lorenz96}.

%
\section{Mathematical background} \label{sec:math}
%

It is well-known \cite{CL83,Carmona} that, assuming sufficient regularity, the desired  closed loop control law is provided by
\begin{equation} \label{eq:optimal control law}
u_t(x) = - R G(x)^{\rm T} \nabla_x v_t(x)
\end{equation}
with the optimal value function $v_t(x)$ satisfying the HJB equation
\begin{equation}  \label{eq:HJB}
-\partial_t v_t = \langle b, \nabla_x v_t \rangle + \frac{1}{2}  \Sigma : 
D_x^2 v_t + c -  \frac{1}{2} \|RG(x)^{\rm T}\nabla_x v_t\|^2_R, \qquad v_T = f.
\end{equation}
Here $A:B = \sum_{i,j} a_{ij}b_{ij}$ and $\langle a,b\rangle = \sum_i a_ib_i$ \cite{Pavliotis2016}. In case of deterministic dynamics; {\it i.e.}~$\Sigma = 0$ in (\ref{eq:FSDE}), the PMP provides the open loop optimal control $U_t = -RG(X_t)^{\rm T}{\rm P}_t$ along the solution ${\rm X}_t$ of (\ref{eq:FSDE}) by solving the Hamiltonian ODEs \cite{holmI}
\begin{subequations} \label{eq:deterministic PMP}
\begin{align}
    \dot{\rm X}_t &= \nabla_P \mathcal{H}_{\rm PMP}({\rm X}_t,{\rm P}_t)= b({\rm X}_t) - G({\rm X}_t)R G({\rm X}_t)^{\rm T} {\rm P}_t,\\
    -\dot{\rm P}_t &= \nabla_X \mathcal{H}_{\rm PMP}({\rm X}_t,{\rm P}_t) =(D_x b({\rm X}_t))^{\rm T}{\rm P}_t - \frac{1}{2} 
    \nabla_x\|R G^{\rm T}({\rm X}_t) {\rm P}_t\|_R^2 + \nabla_x c({\rm X}_t)
    \end{align}
\end{subequations}
in $({\rm X}_t,{\rm P}_t)$ for given Hamiltonian
\begin{equation*}
\mathcal{H}_{\rm PMP}({\rm X}_t,{\rm P}_t) = \langle b({\rm X}_t),{\rm P}_t \rangle + c(X_t)
- \frac{1}{2}\|RG({\rm X}_t)^{\rm T}{\rm P}_t\|_R^2
\end{equation*}
and subject to the boundary conditions ${\rm X}_0 = x_0$ and ${\rm P}_T = \nabla_x f({\rm X}_T)$ \cite{pontryagin}. Formulation (\ref{eq:deterministic PMP}) offers significant computational savings compared to (\ref{eq:HJB}) and has been used, for example, in the context of nonlinear model predictive control (NMPC) \cite{Allgower,MPC}. The PMP for ODEs has been extended to SDEs (\ref{eq:FSDE}) giving rise to the stochastic PMP with (\ref{eq:deterministic PMP}) replaced by a pair of forward-backward SDEs \cite{Carmona,Bensoussan}.

In this paper, instead, we develop a Hamiltonian mean-field formulation in the spirit of the PMP-based Hamiltonian ODE formulation (\ref{eq:deterministic PMP}), which, in turn, can be approximated by interacting particle systems. Since PMP directly targets the gradient of the value function $v_t(x)$; {\it i.e.},~${\rm P}_t = y_t({\rm X}_t)$ with $y_t(x) = \nabla_x v_t(x)$, we note that 
\begin{equation} \label{eq:Pontryagin}
    -\partial_t y_t = \nabla_x \langle b,y_t \rangle - \frac{1}{2}\nabla_x \|RG^{\rm T}y_t\|^2_R +
    \frac{1}{2} \nabla_x (\Sigma : D_x y_t) + \nabla_x c , \qquad
    y_T = \nabla_x f,
\end{equation}
which follows from (\ref{eq:HJB}) assuming sufficient regularity of the value function. 

Substituting the optimal control law (\ref{eq:optimal control law}) into (\ref{eq:FSDE}), the law $\rho_t$ of the states ${\rm X}_t$ satisfies the Fokker--Planck equation
\begin{equation} \label{eq:FPE}
\partial \rho_t = - \nabla_x \cdot \left(\rho_t (b + Gu_t) \right) + \frac{1}{2}
\nabla_x \cdot (\Sigma \nabla_x \rho_t), \qquad \rho_0 = p_0,
\end{equation}
for given initial density $p_0$ \cite{Pavliotis2016}.

We assume throughout this paper that the diffusion matrix $\Sigma$ is position independent. However, the proposed PMP extends to position dependent $\Sigma(x)$. We also allow for $\Sigma =0$, in which case our formulation provides the optimal closed loop control law instead of the open loop solution provided by (\ref{eq:deterministic PMP}). Adding diffusion via $\Sigma = \epsilon I$ with $\epsilon \to 0$ sufficiently small allows one to numerically approximate viscosity solutions of (\ref{eq:HJB}) \cite{CL83}. We will therefore make the following assumptions.

\begin{assumption} \label{assumption}
It is assumed that the HJB equation (\ref{eq:HJB}) possesses  
smooth solutions $v_t$ for given terminal cost $f$. It is also assumed that the associated Fokker--Planck equation (\ref{eq:FPE}) has smooth solutions $\rho_t$ for given initial density $p_0$. It is furthermore assumed that $\rho_t(x)>0$ for all $t\ge0$ and $x\in \mathbb{R}^{d_x}$.
\end{assumption}

%
\section{McKean--Pontryagin Minimum Principle} \label{sec:Hamiltonian}
%

In this section, we develop a novel mean-field approach to the PMP starting from an Eulerian perspective in the spirit of ideal fluid dynamics \cite{Salmon88}. Specifically, we combine the HJB equation (\ref{eq:HJB}) and a Liouville equation for the flow of densities $\rho_t(x)$ under a system of Hamiltonian PDEs. In order to do so, we  assume that the initial states ${\rm X}_0 = x_0$ in (\ref{eq:FSDE}) are no longer fixed but instead distributed according to a given distribution $p_0(x)$. A natural choice for $p_0$ is a Gaussian with mean $x_0$ and a covariance matrix ${\rm C}^{xx}_0$ chosen appropriately. 

Following \cite{Salmon88}, we then transform the Hamiltonian PDEs into a Lagrangian framework in terms of evolving positions $X_t \in \mathbb{R}^{d_x}$ and co-states $P_t \in \mathbb{R}^{d_x}$. These equations are Hamiltonian and provide the novel McKean--Pontryagin formulation. Although we focus on the finite horizon optimal control problem (\ref{eq:costJ}) in this section; an extension to infinite horizon optimal control problems can be found in Section \ref{sec:discounted}. 

%
\subsection{Eulerian formulation}
%

In this subsection, we develop a Hamiltonian formulation of the finite horizon optimal control problem giving rise to evolution equations in the density $\rho_t$ and value function $v_t$, which agree with (\ref{eq:FPE}) and (\ref{eq:HJB}), respectively. 

The first step is to formally decouple the diffusion term in (\ref{eq:HJB}) through an additional function $\phi_t$ and to explicitly reintroduce the closed loop control law $u_t(x)$ into the HJB equation (\ref{eq:HJB}), {\it i.e.},
\begin{subequations}
     \label{eq:HJB2}
     \begin{align}
-\partial_t v_t &= \langle b+Gu_t, \nabla_x v_t\rangle + \frac{1}{2} \Sigma : D_x^2 
\phi_t + c +  \frac{1}{2} \|u_t\|^2_R, \\
0&= \nabla_x \phi_t- \nabla_x v_t,\\
0&= R^{-1}u_t + G^{\rm T} \nabla_x v_t.
\end{align}
\end{subequations}
The second step is to introduce a Lagrange multiplier $\upbeta_t(x)$ that enforces (\ref{eq:HJB2}b) and to define 
the Hamiltonian
\begin{align} \nonumber
    \mathcal{H}(v_t,\rho_t,u_t,\upbeta_t,\phi_t) =&\, \int_{\mathbb{R}^{d_x}}
    \left\{ \bigl< b(x)+G(x)u_t(x),\nabla_x v_t(x)\bigr> + c(x) + \frac{1}{2}\|u_t(x)\|^2_R  \right.\\
    &\,\left.
    +\, \bigl< \upbeta_t(x),\nabla_x (\phi_t(x)-v_t(x))\bigl> + \frac{1}{2} \Sigma : D_x^2\phi_t(x)\right\} \rho_t(x){\rm d}x. \label{eq:Hamiltonian_function_Euler}
\end{align}
We now state the resulting constrained Hamiltonian equations of motion, which arise from (\ref{eq:Hamiltonian_function_Euler}) by taking variational derivatives with respect to its arguments. See  \cite{evans} for the definition of the variational derivative $\updelta \mathcal{L}/\updelta \mu$ of a functional $\mathcal{L}(\mu)$, which depends on a function $\mu$.

\begin{theorem} Let Assumption \ref{assumption} hold. Given the Hamiltonian (\ref{eq:Hamiltonian_function_Euler}), the associated constrained Hamiltonian equations of motion are given by
\begin{subequations} \label{eq:Hamiltonian_Euler}
    \begin{align}
        \partial_t \rho_t &= \frac{\updelta \mathcal{H}}{\updelta v_t} = -\nabla_x \cdot \left( \rho_t (
        b + Gu_t - \upbeta_t)\right),\\
        -\partial_t v_t &= \frac{\updelta \mathcal{H}}{\updelta \rho_t} = \langle b + Gu_t,\nabla_x v_t\rangle + \frac{1}{2} \Sigma: D_x^2 \phi_t +c + \frac{1}{2}\|u_t\|_R^2
        ,\\
        0 &= \frac{\updelta \mathcal{H}}{\updelta \upbeta_t} = \rho_t \left( \nabla_x \phi_t-\nabla_x v_t\right),\\
        0 &= \frac{\updelta \mathcal{H}}{\updelta \phi_t} = -\nabla_x\cdot\left(
        \rho_t \left( \upbeta_t - \frac{1}{2}\Sigma \nabla_x \log\rho_t\right)\right),\\
        0&= \frac{\updelta \mathcal{H}}{\updelta u_t}= \rho_t \left( G^{\rm T}\nabla_x v_t + R^{-1}u_t\right).
    \end{align}
\end{subequations}
The time evolutions of $\rho_t$ and $v_t$ agree with those given by the Fokker--Planck equation (\ref{eq:FPE}) and the HJB equation (\ref{eq:HJB}), respectively.
The Hamiltonian (\ref{eq:Hamiltonian_function_Euler}) is conserved along solutions of (\ref{eq:Hamiltonian_Euler}); {\it i.e.},
\begin{equation*}
\frac{{\rm d}}{{\rm d}t}\mathcal{H}(v_t,\rho_t,u_t,\upbeta_t,\phi_t) = 0.
\end{equation*}

\end{theorem}

\begin{proof}
    Substituting (\ref{eq:Hamiltonian_Euler}d) into (\ref{eq:Hamiltonian_Euler}a) leads to the Fokker--Planck equation (\ref{eq:FPE}). Similarly, substituting (\ref{eq:Hamiltonian_Euler}c) and (\ref{eq:Hamiltonian_Euler}e) into (\ref{eq:Hamiltonian_Euler}b) leads to the HJB equation (\ref{eq:HJB}). We also have
    \begin{align*}
    &\frac{{\rm d}}{{\rm d}t}\mathcal{H}(v_t,\rho_t,u_t,\upbeta_t,\phi_t) = \\
    &\qquad \int_{\mathbb{R}^{d_x}} \left\{ \frac{\updelta \mathcal{H}}{\updelta v_t}\partial_t v_t + \frac{\updelta \mathcal{H}}{\updelta \rho_t}\partial_t \rho_t + 
    \Bigl< \frac{\updelta \mathcal{H}}{\updelta \upbeta_t},\partial_t \upbeta_t\Bigr> \,+\,
    \frac{\updelta \mathcal{H}}{\updelta \phi_t}\partial_t \phi_t +
    \Bigl< \frac{\updelta \mathcal{H}}{\updelta u_t},\partial_t u_t \Bigr>\right\} {\rm d}x=0,
    \end{align*}
    which implies conservation of the Hamiltonian (\ref{eq:Hamiltonian_function_Euler}) along solutions of (\ref{eq:Hamiltonian_Euler}).
\end{proof}

\begin{remark} \label{rem0}
    It should be noted that $\phi_t$ and $\upbeta_t$ are not uniquely determined by constraints (\ref{eq:Hamiltonian_Euler}c)-(\ref{eq:Hamiltonian_Euler}d). Indeed, any time-dependent constant $c_t$ can be added to $v_t$; {\it i.e.}
    \begin{equation*}
    \phi_t(x) = v_t(x) + c_t,
    \end{equation*}
    and any time-dependent function $\upbeta'_t(x)$ with $\nabla_x \cdot (\rho_t\upbeta_t') = 0$ to $\frac{1}{2}\Sigma \nabla_x \log \rho_t$; {\it i.e.}
    \begin{equation} \label{eq:gauge beta freedom}
        \upbeta_t(x) = \frac{1}{2}\nabla_x \log \rho_t(x) + \upbeta'_t(x),
    \end{equation}
    while not changing the evolution equations in $\rho_t$ and $v_t$.
    This fact constitutes a gauge freedom of the proposed formulation. A more general type of gauge invariance will be exploited in Subsection \ref{sec:gauge}, where $\upbeta_t$ will be allowed to modify the time evolution of the density $\rho_t$ while simultaneously preserving the time evolution in the value function $v_t$.
\end{remark}

Sometimes, it is useful to have the corresponding action functional available, which also encodes the given boundary conditions $\rho_0 = p_0$ and $v_T = f$.

\begin{lemma}
The action functional leading to (\ref{eq:Hamiltonian_Euler}) is given by
\begin{align} \nonumber
    \mathcal{S}(v,\rho,u,\upbeta,\phi) =& \int_0^T \left\{
    \int_{\mathbb{R}^{d_x}}  \rho_t \partial_tv_t \,{\rm d}x + 
    \mathcal{H}(v_t,\rho_t,u_t,\upbeta_t,\phi_t)\right\}{\rm d}t \\
    &\,+\,
    \int_{\mathbb{R}^{d_x}} \left\{ (f-v_T)\rho_T + v_0 p_0\right\} {\rm d}x. \label{eq:action Euler}
\end{align}
\end{lemma}

\begin{proof} Equations (\ref{eq:Hamiltonian_Euler}) follow from (\ref{eq:action Euler}) using integration by parts; {\it i.e.},
\begin{equation}
    \int_0^T v_t\partial_t \rho_t + \int_0^T 
    \rho_t \partial_t v_t = v_T \rho_T - v_0\rho_0.
\end{equation}
Variations with respect to $\rho_T$ and $v_0$ lead to the desired boundary conditions.
\end{proof}

\begin{remark} \label{rem1}
A related variational saddle point formulation  has been proposed and exploited numerically in \cite{MLLO24} for solving Hamilton--Jacobi partial differential equations. In the context of our problem setting with $G \equiv I$, $b \equiv 0$, and $\Sigma = 2\epsilon I$, the formulation stated in \cite{MLLO24} leads to the action functional
\begin{subequations} 
\begin{align*}
\widetilde{\mathcal{S}}(\phi,\rho,u) =& \int_0^T  
\int_{\mathbb{R}^{d_x}}  \rho_t \left( 
\partial_t \phi_t + \epsilon D_x^2 \phi_t +
\langle u_t,\nabla_x \phi_t \rangle + c + \frac{1}{2}\|u_t\|_R^2
\right\}{\rm d}x  {\rm d}t \\ &\,+\,
\int_{\mathbb{R}^{d_x}} 
\left\{ (f-\phi_T)\rho_T + \gamma \phi_0\right\}\,{\rm d}x.
\end{align*}
\end{subequations}
Here, $\gamma>0$ is a hyperparameter and the optimal control law $u_t(x) = -\nabla_x \phi_t$ is obtained by minimizing $\widetilde{\mathcal{S}}$ over $u_t$. The implied evolution equations in $\rho_t$ and $\phi_t$  agree with (\ref{eq:Hamiltonian_Euler}a)-(\ref{eq:Hamiltonian_Euler}b) for $\upbeta_t = \frac{1}{2} \Sigma \nabla_x \log \rho_t$ and $v_t = \phi_t$. However, the initial condition for $\rho_t$ is $\rho_0 = \gamma$.
\end{remark}

\noindent
In the following subsection, we demonstrate how our Hamiltonian formulation (\ref{eq:Hamiltonian_Euler}) facilitates the transition to the Lagrangian perspective and provides enhanced flexibility due to the underlying gauge invariance, as further discussed in Section \ref{sec:gauge}.

%
\subsection{Lagrangian formulation}
%

We now convert (\ref{eq:Hamiltonian_Euler}) into the corresponding Lagrangian formulation \cite{Salmon88}. Let $X_t(a)$ denote the solution of the mean-field ODE
\begin{equation} \label{eq:ODE}
\dot{X}_t(a) = b(X_t(a)) + G(X_t(a))u_t(X_t(a))
- \upbeta_t(X_t(a))
\end{equation}
for given control law $u_t(x)$ and $\upbeta_t(x)$ of the form (\ref{eq:gauge beta freedom}). Following \cite{Salmon88}, the labels $a \in \Omega := [0,1]^{d_x}$ are chosen so that
\begin{equation}
    p_0(X_0(a))\,|D_a X_0(a)| =1.
\end{equation}
Here, $|D_a X_0(a)|$ denotes the absolute value of the determinant of the Jacobian matrix $D_a X_0(a) \in \mathbb{R}^{d_x \times d_x}$. We denote the thus defined map from labels $a\in \Omega$ to initial conditions $X_0(a)\in \mathbb{R}^{d_x}$ by $\Psi:\Omega\to \mathbb{R}^{d_x}$; {\it i.e.},
\begin{equation} \label{eq:ICmap}
    X_0(a) = \Psi(a).
\end{equation}

Since the labels $a \in \Omega$ are independent variables and fixed along trajectories $X_t(a)$, we have
\begin{equation}
\rho_t(X_t(a)) \,|D_a X_t(a)| = 1
\end{equation}
for all $t >0$ \cite{Salmon88}, where $\rho_t$ satisfies the Liouville equation (\ref{eq:Hamiltonian_Euler}a). It follows that
\begin{equation} \label{eq:integral conversion}
    \int_{\mathbb{R}^{d_x}} g(x) \,\rho_t(x)\,{\rm d}x = 
    \int_\Omega g(X_t(a)) \,{\rm d}a
    = \int_\Omega g_t(a) \,{\rm d}a
\end{equation}
for any continuous function $g:\mathbb{R}^{d_x}\to \mathbb{R}$ and
$g_t(a) = g(X_t(a))$. We refer to $\rho_t$ as the density of $X_t$ and write $X_t \sim \rho_t$. The particle labeling (\ref{eq:ICmap}) implies $X_0 \sim p_0$ at initial time.

We now reformulate the Hamiltonian (\ref{eq:Hamiltonian_function_Euler}) in Lagrangian coordinates giving rise to
\begin{align} \nonumber
\mathrm{H}(P_t,X_t,U_t,\beta_t,\phi_t) =&
\int_\Omega \biggl\{ \bigl< 
b(X_t(a)) + G(X_t(a))U_t(a),P_t(a) \bigr> \biggr. \\ \nonumber
&\,+\, c(X_t(a)) + \frac{1}{2}\|U_t(a)\|_R^2  \\
&\,\left. + \bigl< \beta_t(a), \nabla_x \phi_t(X_t(a))-P_t(a)\bigr>
+ \frac{1}{2} \Sigma : D_x^2 \phi_t(X_t(a))
\right\}{\rm d}a. \label{eq:Hamiltonian_function_Lagrange}
\end{align}
Here we have introduced co-states $P_t(a)$, controls $U_t(a)$, and Lagrange multipliers $\beta_t(a)$. 

\begin{theorem} \label{theorem2}
Let Assumption \ref{assumption} hold.
Given the Hamiltonian (\ref{eq:Hamiltonian_function_Lagrange}), the associated constrained Hamiltonian equations of motion are given by
\begin{subequations} \label{eq:Hamiltonian_Lagrange}
    \begin{align}
    \dot{X}_t &= \frac{\updelta \mathrm{H}}{\updelta P_t} =
    b(X_t) + G(X_t)U_t - \beta_t ,\\ \nonumber
        -\dot{P}_t &= \frac{\updelta \mathrm{H}}{\updelta X_t} = \left(D_x(b(X_t)+G(X_t)U_t)\right)^{\rm T}P_t + \nabla_x c(X_t) \\
        & \quad \qquad \quad \,\,\,\,+ \,\frac{1}{2}\nabla_x (\Sigma:D_x^2 \phi_t)(X_t) + D_x^2 \phi_t(X_t)\beta_t,\\
        0 &= \frac{\updelta \mathrm{H}}{\updelta \beta_t} =  \nabla_x \phi_t(X_t) -P_t,\\
        0 &= \frac{\updelta \mathrm{H}}{\updelta \phi_t} = -\nabla_x\cdot\left(
        \rho_t \left(\beta_t - \frac{1}{2}\Sigma \nabla_x \log\rho_t(X_t)\right)\right),\\
        0&= \frac{\updelta \mathrm{H}}{\updelta U_t}=  G(X_t)^{\rm T}P_t + R^{-1}U_t.
    \end{align}
\end{subequations}
The Hamiltonian (\ref{eq:Hamiltonian_function_Lagrange}) is conserved along solutions of (\ref{eq:Hamiltonian_Lagrange}); {\it i.e.},
\begin{equation*}
\frac{{\rm d}}{{\rm d}t}\mathrm{H}(P_t,X_t,U_t,\beta_t,\phi_t) = 0.
\end{equation*}

Given common boundary conditions
\begin{equation*}
    \rho_0 = p_0, \quad
    X_0 \sim p_0, \quad v_T = f, \quad P_t = \nabla_x f(X_t), 
\end{equation*}
the solution $(v_t,\rho_t,u_t,\upbeta_t,\phi_t)$ of (\ref{eq:Hamiltonian_Euler}) and the solution $(X_t,P_t,U_t,\beta_t,\phi_t)$ of (\ref{eq:Hamiltonian_Lagrange}) 
satisfy
\begin{equation} \label{eq:correspondence}
X_t \sim \rho_t, \quad P_t = \nabla_x v_t(X_t), \quad U_t = u_t(X_t)
\end{equation}
for $t\in [0,T]$. 
Furthermore, the Hamiltonian functionals (\ref{eq:Hamiltonian_function_Euler}) and (\ref{eq:Hamiltonian_function_Lagrange}) agree; {\it i.e.},
\begin{equation*}
    \mathrm{H}(P_t,X_t,U_t,\beta_t,\phi_t) = \mathcal{H}(v_t,\rho_t,u_t,\upbeta_t,\phi_t)
\end{equation*}
\end{theorem}

\begin{proof}
We have
    \begin{align*}
    &\frac{{\rm d}}{{\rm d}t}\mathrm{H}(P_t,X_t,U_t,\beta_t,\phi_t) =\\
    &\qquad \int_\Omega \left\{ \Bigl< \frac{\updelta \mathrm{H}}{\updelta P_t},\dot{P}_t \Bigr> + \Bigl< \frac{\updelta \mathrm{H}}{\updelta X_t},\dot{X}_t \Bigr>  + 
    \Bigl< \frac{\updelta \mathrm{H}}{\updelta \beta_t},\dot{\beta}_t \Bigr> \,+\,
    \frac{\updelta \mathrm{H}}{\updelta \phi_t}\partial_t \phi_t +
    \Bigl< \frac{\updelta \mathrm{H}}{\updelta U_t}\dot{U}_t \Bigr> \right\}{\rm d}a 
    =0,
    \end{align*}
    which implies conservation of the Hamiltonian (\ref{eq:Hamiltonian_function_Lagrange}) along solutions of (\ref{eq:Hamiltonian_Lagrange}).
    
Since the mean-field ODE (\ref{eq:Hamiltonian_Lagrange}a) has (\ref{eq:Hamiltonian_Euler}a) as its associated Liouville equation, the mutual relationships expressed in (\ref{eq:correspondence}) all follow from a direct comparison of (\ref{eq:Hamiltonian_Euler}) and (\ref{eq:Hamiltonian_Lagrange}) once it has been established from (\ref{eq:Hamiltonian_Lagrange}) that $\nabla_x \phi_t = \nabla_x v_t = y_t$ where $y_t$ satisfies (\ref{eq:Pontryagin}). Therefore, let us derive the evolution equation for $\nabla_x \phi_t(x)$ implied by (\ref{eq:Hamiltonian_Lagrange}):
\begin{align}\label{eq:proof1} \nonumber
        -\partial_t \nabla_x \phi_t(X_t) =& \,\,D_x^2 \phi_t(X_t)\dot{X}_t -\dot{P}_t\\ \nonumber
        =& \,\,D_x^2 \phi_t(X_t)\left( b(X_t) + G(X_t)U_t - \beta_t\right) \\ \nonumber
        &+ \left(D_x (b(X_t)+G(X_t)U_t)\right)^{\rm T}P_t + \nabla_x c(X_t)
        \\ \nonumber
        &+ \frac{1}{2} \nabla_x  (\Sigma : D_x^2\phi_t)(X_t)  
        + D_x^2 \phi_t(X_t)\beta_t  \\ \nonumber
        =& \,\,\nabla_x \langle b(X_t),\nabla_x \phi_t(X_t)\rangle  -\frac{1}{2} \nabla_x \|R G(X_t)^{\rm T}\nabla_x \phi_t(X_t)\|^2_R 
        \\
        &+ \nabla_x c(X_t) + \frac{1}{2} \nabla_x  (\Sigma : D_x^2 \phi_t)(X_t) ,  
\end{align}
where we have used (\ref{eq:Hamiltonian_Lagrange}c) and 
\begin{equation*} 
U_t = -RG(X_t)^{\rm T}\nabla_x \phi_t(X_t).
\end{equation*}
Also note that
\begin{align*}
    \frac{1}{2} \nabla_x \|R G(X_t)^{\rm T}\nabla_x \phi_t(X_t)\|^2_R &=
    -\nabla_x \langle G(X_t)U_t,\nabla_x \phi_t(X_t)\rangle \\
    &= 
    -\left(D_x G(X_t)U_t \right)^{\rm T}\nabla_x \phi_t(X_t) - 
    D_x^2 \phi_t(X_t)G(X_t)U_t.
\end{align*}
Hence, it follows that $y_t :=\nabla_x \phi_t$ satisfies (\ref{eq:Pontryagin}). Furthermore, $\nabla_x \phi_T = \nabla_x f(X_T)$. 

We already emphasize at this point that all terms involving $\beta_t$ cancel out, and therefore the calculations in (\ref{eq:proof1}) are valid for any choice of $\beta_t$ as long as the implied flow map $X_t(a)$ remains smooth and invertible. 

The equivalence of the Hamiltonian functions $\mathcal{H}$ and $\mathrm{H}$ along solutions follows from (\ref{eq:integral conversion}) and (\ref{eq:correspondence}). 
\end{proof}

\begin{remark}
    One may view the co-states $P_t$ as momenta in ideal fluid dynamics in which case the relationship $P_t(a) = \nabla_x v_t(X_t(a))$ implies that the fluid is irrotational \cite{Salmon88}.
\end{remark}

Similarly to the action principle (\ref{eq:action Euler}) in the Eulerian setting, there is a corresponding action principle for the constrained Hamiltonian formulation (\ref{eq:Hamiltonian_Lagrange}). 

\begin{lemma}
The action functional leading to (\ref{eq:Hamiltonian_Lagrange}) is given by
\begin{align} \nonumber
    \mathrm{S}(P,X,U,\beta,\phi) =& \int_0^T \left\{
    \int_\Omega  \bigl< \dot{P}_t(a),X_t(a) \bigr>\,
    {\rm d}a + \mathrm{H}(P_t,X_t,U_t,\beta_t,\phi_t)\right\}{\rm d}t \\
    &\,+\,
    \int_\Omega \left(f(X_T(a)) - \bigl< P_T(a),X_T(a)\bigr>
    -\bigl< P_0(a),\Psi(a) \bigr> \right)  {\rm d}a. \label{eq:action Lagrange}
\end{align}
\end{lemma}

\begin{proof}
    The Hamiltonian equations of motion (\ref{eq:Hamiltonian_Lagrange}) follow from the action (\ref{eq:action Lagrange}) and integration by parts
    \begin{equation*}
    \int_0^T \bigl< P_t(a),\dot{X}_t(a) \bigr> \,{\rm d}t
    + \int_0^T \bigl< \dot{P}_t(a), X_t(a) \bigr> \,{\rm d}t
    = \bigl< P_T(a),X_T(a)\bigr> - \bigl< P_0(a),X_0(a) \bigr>
    \end{equation*}
    In particular, we obtain the boundary conditions $P_T(a) = \nabla_x f(X_T(a))$ and (\ref{eq:ICmap}).
\end{proof}

Upon introducing the generator
\begin{equation*} 
\mathcal{L}_{\upbeta_t} g = \bigl< \upbeta_t , \nabla_x g \bigr>
+ \frac{1}{2} \Sigma : D_x^2 g
\end{equation*}
for suitable functions $g(x)$, we have the following
\begin{equation} \label{eq:generator rewrite}
\int_\Omega \left\{
\bigl< \beta_t,\nabla_x \phi_t(X_t)\bigr> + \frac{1}{2}
\Sigma : D_x^2 \phi_t(X_t)\right\} {\rm d}a
= \int_{\mathbb{R}^{d_x}} (\mathcal{L}_{\upbeta_t}\phi_t) \,\rho_t\,{\rm d}x,
\end{equation}
which allows us to rewrite the equations of motion (\ref{eq:Hamiltonian_Lagrange}a)-(\ref{eq:Hamiltonian_Lagrange}b). More specifically, we introduce the pullback of $\mathcal{L}_{{\upbeta}_t}$ to functions over label space $\Omega$ via the flow map $X_t(a)$, which we denote by $\widehat{\mathcal{L}}_{{\beta}_t}$; {\it i.e.}, 
\begin{equation} \label{eq:definition pulback operator}
\widehat{\mathcal{L}}_{{\beta}_t} \widehat{g}_t = \mathcal{L}_{{\upbeta}_t}g (X_t)
\end{equation}
with $\widehat{g}_t(a) = g(X_t(a))$ and $\beta_t(a) = \upbeta_t(X_t(a))$. Finally, it should be noted that 
\begin{equation*}
    \frac{1}{2}\nabla_x (\Sigma:D_x^2 \phi_t)(X_t) + D_x^2 \phi_t(X_t)\beta_t = \mathcal{L}_{{\upbeta}_t} \nabla_x \phi_t(X_t) =  \widehat{\mathcal{L}}_{{\beta}_t} P_t
\end{equation*}
since $P_t(a) = \nabla_x \phi_t(X_t(a))$.

With these notations in place, we can reformulate the evolution equations (\ref{eq:Hamiltonian_Lagrange}a)-(\ref{eq:Hamiltonian_Lagrange}b) as
\begin{subequations} \label{eq:Hamiltonian ODE}
    \begin{align}
        \dot{X}_t =&  \,\,b(X_t) - G(X_t)R G(X_t)^{\rm T}P_t - \widehat{\mathcal{L}}_{\beta_t}X_t,\\ 
        -\dot{P}_t 
        =& \,\left(D_x b(X_t)\right)^{\rm T}P_t 
        - \frac{1}{2} \nabla_x 
        \|RG(X_t)^{\rm T}P_t\|^2_R + \nabla_x c(X_t) + \widehat{\mathcal{L}}_{\beta_t} P_t,
    \end{align}
\end{subequations}
which are to be compared to the classical PMP formulation (\ref{eq:deterministic PMP}). The only new terms arise through the generator $\widehat{\mathcal{L}}_{\beta_t}$. 

In line with (\ref{eq:gauge beta freedom}), a natural choice for the function $\upbeta_t$ is to set $\upbeta'_t(x)=0$, which results in
\begin{equation} \label{eq:beta natural}
    \widetilde{\upbeta}_t(x) = \frac{1}{2}\Sigma \nabla_x \log \rho_t(x),
\end{equation}
where $\rho_t(x)$ is the density of $X_t$ as defined by the Liouville equation (\ref{eq:Hamiltonian_Euler}a). This choice leads to
\begin{equation} \label{eq:natural generator}
    \mathcal{L}_{\widetilde{\upbeta}_t}g =  \frac{1}{2}\bigl< \Sigma \nabla_x \log \rho_t , \nabla_x g \bigr> + \frac{1}{2} \Sigma : D_x^2 g,
\end{equation}
which is the generator of a reversible diffusion process with invariant measure $\rho_t$ \cite{Pavliotis2016}.  Since
\begin{equation*}
\int_\Omega  \widehat{\mathcal{L}}_{\widetilde{\beta}_t} \widehat{\phi}_t(a) \,{\rm d}a
= \int_{\mathbb{R}^{d_x}} (\mathcal{L}_{\widetilde\upbeta_t} \phi_t)(x) \,\rho_t(x)\,{\rm d}x = 0,
\end{equation*}
$\widehat{\phi}_t(a) = \phi_t(X_t(a))$, and (\ref{eq:generator rewrite}), the Hamiltonian (\ref{eq:Hamiltonian_function_Lagrange}) simplifies to
\begin{align} \nonumber
\widetilde{\mathrm{H}}(P_t,X_t,U_t) =&
\int_\Omega \biggl\{ \bigl< 
b(X_t(a)) + G(X_t(a))U_t(a),P_t(a) \bigr> \biggr. \\ 
&\,+\, c(X_t(a)) + \frac{1}{2}\|U_t(a)\|_R^2  - \bigl<\widehat{\mathcal{L}}_{\widetilde{\beta}_t}X_t(a),P_t(a)\bigr>
\biggr\}\,{\rm d}a \label{eq:Hamiltonian_function_Lagrange_generator_simple}
\end{align}
under the choice $\beta_t = \widetilde{\beta}_t$. 

\begin{remark}
 This paper introduces a deterministic analog of the forward-back- ward SDE formulation of stochastic optimal control in terms of PMP. Alternatively, \cite{ZOL24} provides a corresponding reformulation in terms of  mean-field differential equations in $(X_t(a),Y_t(a)) \in \mathbb{R}^{d_x}\times \mathbb{R}$ and the value function $\phi_t$ is obtained from
\begin{equation}
\phi_t(X_t(a)) = Y_t(a);
\end{equation}
{\it i.e.}, the backward differential equation of $Y_t$ states a Lagrangian evolution of the value function instead of its gradient. This is in contrast to our formulation which relies on co-states $P_t(a) \in \mathbb{R}^{d_x}$, as in (\ref{eq:Hamiltonian_Lagrange}c) and leads to a Hamiltonian formulation in $(X_t,P_t)$. The forward evolution equations of $X_t$ agree for both formulations if $\beta_t = \widetilde{\beta}_t$
in (\ref{eq:Hamiltonian_Lagrange}d).
\end{remark}

We close this subsection by reemphasizing the fact that the evolution equations (\ref{eq:Hamiltonian ODE}) are of mean-field type through the appearance of the generator $\mathcal{L}_{\upbeta_t}$ and the function $\phi_t$, which is implicitly defined through the constraint (\ref{eq:Hamiltonian_Lagrange}c), although we are solving a classical stochastic optimal control problem. Extensions to mean-field control problems are feasible; but beyond the scope of this contribution. See, \cite{R25c} for a first application in the context of optimal control of partially observed diffusion processes. 

%
\subsection{Gauge freedom} \label{sec:gauge}
%

In this subsection, we further explore the freedom to choose $\beta_t(a) = \upbeta_t(X_t(a))$ in the context of finite horizon optimal control problems. 
Recall the gauge freedom discussed in Remark \ref{rem0}. Here we extend the gauge freedom by only requesting that the time evolution of $\phi_t$ and the optimal control law
\begin{equation*}
u_t(x) = -R G(x)^{\rm T}\nabla_x \phi_t(x).
\end{equation*}
remain unaffected by the choice of $\beta_t$, while the time evolution of the particle distribution $\rho_t$ can be modified.

Indeed, consider the reduced set of constrained Hamiltonian equations of motion
\begin{subequations} \label{eq:Hamiltonian_Lagrange_reduced}
    \begin{align}
    \dot{X}_t &= 
    b(X_t) + G(X_t)U_t - \beta_t ,\\ \nonumber
        -\dot{P}_t &= \left(D_x(b(X_t)+G(X_t)U_t)\right)^{\rm T}P_t + \nabla_x c(X_t) \\
        & \quad   \,\,+ \,\frac{1}{2}\nabla_x (\Sigma:D_x^2 \phi_t)(X_t) + D_x^2 \phi_t(X_t)\beta_t,\\
        0 &=   \nabla_x \phi_t(X_t) -P_t,\\
        0&=   G(X_t)^{\rm T}P_t + R^{-1}U_t.
    \end{align}
\end{subequations}
Following the proof of Theorem \ref{theorem2} and the calculations (\ref{eq:proof1}) in particular, we find that $\nabla_x \phi_t$ satisfies (\ref{eq:Pontryagin}) for any choice of $\beta_t$ as long as the induced flow map $X_t(a)$ remains smooth and invertible. Hence, the desired optimal control is still provided by (\ref{eq:HJB2}c). Of course, the Liouville equation (\ref{eq:Hamiltonian_Euler}a) and the time evolution of $X_t$ are altered. We summarize this finding in the following theorem.

\begin{theorem} \label{theorem1}
Let assumption \ref{assumption} hold. For any choice of the variable $\beta_t$ for which the induced flow map $X_t(a)$ remains smooth and invertible, the time evolution of $\phi_t(x)$ resulting from (\ref{eq:Hamiltonian_Lagrange_reduced}) with terminal condition $P_t = \nabla_x f(X_t)$ satisfies
\begin{equation*}
    \nabla_x \phi_t(x) = \nabla_x v_t(x),
\end{equation*}
where $v_t(x)$ is the solution of the HJB equation (\ref{eq:HJB}).
\end{theorem}

\begin{proof}
We have already demonstrated as part of the proof of Theorem \ref{theorem2} that $\nabla_x \phi_t$ satisfies the HJB equation (\ref{eq:Pontryagin}) and that the contributions of $\beta_t$ on the time evolution of $X_t$ and $P_t$ cancel each other out in the time-differentiated constraint (\ref{eq:Hamiltonian_Lagrange_reduced}c):
\begin{equation}
    0 = \partial_t \nabla_x \phi_t(X_t) + D_x^2 \phi_t(X_t)\dot{X}_t - \dot{P}_t;
\end{equation} 
compare (\ref{eq:proof1}).
\end{proof}

Consider, for example, 
\begin{equation} \label{eq:beta decouple}
\beta_t = -G(X_t) R G(X_t)^{\rm T} P_t,
\end{equation}
which eliminates the control from the evolution equation in $X_t$ since
\begin{equation*}
G(X_t) U_t = -G(X_t) RG(X_t)^{\rm T} P_t = \beta_t.
\end{equation*}
Hence, the boundary value problem decouples and a single forward and reverse time simulation of 
\begin{subequations} 
    \begin{align*}
        \dot{X}_t =& \,\,b(X_t) ,\\
        -\dot{P}_t =& \,\left(D_x b(X_t) \right)^{\rm T}P_t 
        -\frac{1}{2} \nabla_x \|RG(X_t)^{\rm T} P_t\|^2_R
        \\ &+ \nabla_x c(X_t) + \frac{1}{2} \nabla_x  (\Sigma : D_x^2 \phi_t)(X_t)
        + D_x^2 \phi_t(X_t)\beta_t,\\
        0 =&\,\,\nabla_x \phi_t(X_t)-P_t
    \end{align*}
\end{subequations}
with initial condition (\ref{eq:ICmap}) and terminal condition $P_T(a) = \nabla_x f(X_T(a))$ solves the finite horizon stochastic optimal control problem. We note that such a decoupling is also possible when $\Sigma$ is singular or even zero. This is in contrast to classical and stochastic PMP \cite{Carmona}.

A refined choice is provided by
\begin{equation*}
    \beta_t = \frac{1}{2}\Sigma \nabla_x \log \rho_t(X_t) -
    G(X_t)RG(X_t)^{\rm T} P_t - G(X_t) u_t^{\rm ref}(X_t)  ,
\end{equation*}
where $u_t^{\rm ref}(x)$ is a known reference control law leading to
\begin{subequations} \label{eq:NMPC formulation}
    \begin{align}
        \dot{X}_t =& \,\,b(X_t) + G(X_t) u_t^{\rm ref}(X_t) -
        \widehat{\mathcal{L}}_{\widetilde{\beta}_t}(X_t),\\ \nonumber
        -\dot{P}_t =& \,\left(D_x b(X_t) \right)^{\rm T}P_t 
        -\frac{1}{2} \nabla_x \|RG(X_t)^{\rm T} P_t\|^2_R + \nabla_x c(X_t)
        \\ & \,+\, \widehat{\mathcal{L}}_{\widetilde{\beta}_t}P_t
        - D_x^2 \phi_t(X_t)\,G(X_t)\left(RG(X_t)^{\rm T} P_t +  u_t^{\rm ref}(X_t)\right),
        \\
        0 =&\,\,\nabla_x \phi_t(X_t)-P_t,
    \end{align}
\end{subequations}
where we have used the pulled back generator $\widehat{\mathcal{L}}_{\widetilde{\beta}_t}$; compare (\ref{eq:definition pulback operator}). In the context of NMPC \cite{MPC,Allgower}, such a control could arise from the previous receding horizon control window.

\begin{remark}
Formally, it is also possible to set $\beta_t = -\Sigma^{1/2}\dot{B}_t$ in (\ref{eq:Hamiltonian_Lagrange}). The evolution equations in $X_t$ and $P_t$ then become SDEs, which should be interpreted in Stratonovich form \cite{Pavliotis2016}. More precisely, (\ref{eq:Hamiltonian_Lagrange}a) becomes the SDE (\ref{eq:FSDE}) while (\ref{eq:Hamiltonian_Lagrange}b) becomes
\begin{align*}
-\dot{P}_t &= \frac{\updelta \mathrm{H}}{\updelta X_t} = \left(D_x(b(X_t)+G(X_t)U_t)\right)^{\rm T}P_t + \nabla_x c(X_t) \\
        & \quad \qquad \quad \,\,\,\,+ \,\frac{1}{2}\nabla_x (\Sigma:D_x^2 \phi_t)(X_t) - D_x^2 \phi_t(X_t)\Sigma^{1/2} \circ \dot{B}_t.
\end{align*}
We still have that $\phi_t$, as determined by (\ref{eq:Hamiltonian_Lagrange}c), undergoes a deterministic time evolution in agreement with (\ref{eq:Pontryagin}). This formulation is related to the forward-backward SDE formulation of stochastic optimal control \cite{Carmona}. An exploration of this connection is left to future research.
\end{remark}

%
\section{Discounted cost optimal control problems} \label{sec:discounted}
%

We now consider infinite horizon optimal control problems with discounted cost function (\ref{eq:costJ infinite}). The associated stationary HJB equation is given by
\begin{equation} \label{eq:stationary HJB}
0= -\gamma \hat v + \langle b, \nabla_x \hat v \rangle + \frac{1}{2}  \Sigma 
: D_x^2 \hat v + c -  \frac{1}{2} \|RG(x)^{\rm T}\nabla_x \hat v\|^2_R .
\end{equation}
We assume that $\hat v:\mathbb{R}^{d_x}\to \mathbb{R}$ exists, is unique, and is sufficiently regular. 
Let us denote the associated optimal control law by $\hat u$; {\it i.e.},
\begin{equation*}
    \hat u(x) = -RG(x)^{\rm T}\nabla_x \hat v(x).
\end{equation*}

We furthermore consider the forward HJB equation
\begin{equation}
     \label{eq:HJB2d}
\partial_t v_t =  -\gamma v_t + \langle b , \nabla_x v_t\rangle + \frac{1}{2}  \Sigma : D_x^2 v_t + c -  \frac{1}{2} \|RG^{\rm T}\nabla_x v_t\|^2_R, \qquad v_0 = 0,
\end{equation}
and make the following assumptions.

\begin{assumption} \label{assumption4}
    We assume that the solution $v_t$, $t\ge 0$, of the HJB equation (\ref{eq:HJB2d}) is sufficiently regular and that the control law 
    \begin{equation*}
        u_t(x) = -RG(x)^{\rm T} \nabla_x v_t(x)
    \end{equation*}
    converges to the control law $\hat u(x)$ for $t\to \infty$. Furthermore, we assume that the density $\rho_t$ defined by the associated Fokker--Planck equation (\ref{eq:FPE}) with $\rho_0 = p_0$ converges to a unique invariant distribution denoted $\hat \rho$ as $t\to \infty$. 
\end{assumption}

We adjust our McKean--Pontryagin formulation (\ref{eq:Hamiltonian_Lagrange_reduced})
to the problem of infinite horizon optimal control. The first step is to revert time and add the damping term $-\gamma P_t$ to the evolution equation (\ref{eq:Hamiltonian_Lagrange_reduced}b). This results in constrained evolution equations
\begin{subequations} \label{eq:time-reversed PMP}
\begin{align}
    \dot{X}_t &=
    -b(X_t) - G(X_t)U_t + \beta_t ,\\ \nonumber
        \dot{P}_t &= -\gamma P_t + \left(D_x(b(X_t)+G(X_t)U_t)\right)^{\rm T}P_t + \nabla_x c(X_t) \\
        & \quad   \,\,+ \,\frac{1}{2}\nabla_x (\Sigma:D_x^2 \phi_t)(X_t) + D_x^2 \phi_t(X_t)\beta_t,\\
        0 &=   \nabla_x \phi_t(X_t) -P_t,\\
        0&=  G(X_t)^{\rm T}P_t + R^{-1}U_t.
\end{align}
\end{subequations}
We next choose $\beta_t$ so that the forward evolution of $X_t$ leads to a Liouville equation with $\hat \rho$ as the stationary distribution as $t\to \infty$. This requirement leads to
\begin{equation} \label{eq:time reversed beta}
    \beta_t = 2b(X_t) + 2G(X_t)U_t - \frac{1}{2}\Sigma \nabla_x \log \rho_t(X_t).
\end{equation}
We summarize the resulting equations of motion using the generator $\widehat{\mathcal{L}}_{\widetilde{\beta}_t}$:
\begin{subequations} \label{eq:PMP_forward}
\begin{align}
    \dot{X}_t &=b(X_t) + G(X_t)U_t - \widehat{\mathcal{L}}_{\widetilde{\beta}_t}X_t ,\\ \nonumber
        \dot{P}_t &= -\gamma P_t + \left(D_x(b(X_t)+G(X_t)U_t)\right)^{\rm T}P_t + \nabla_x c(X_t) \\
        & \quad   \,\,+ \,\widehat{\mathcal{L}}_{\widetilde{\beta}_t}P_t + 2D_x^2 \phi_t(X_t)\dot{X}_t,\\
        0 &=   \nabla_x \phi_t(X_t) -P_t,\\
        0&=  G(X_t)^{\rm T}P_t + R^{-1}U_t.
\end{align}
\end{subequations}

\begin{remark}
Assume that initial conditions for (\ref{eq:PMP_forward}) satisfy $X_0\sim \hat \rho$ and $P_0(a) = \nabla_x \hat v(X_0(a))$, then it follows that the resulting $\phi_t$ satisfies the stationary HJB equation (\ref{eq:stationary HJB}) possibly up to an irrelevant constant $c$; {\it i.e.}, $\phi_t = \hat v + c$, for all $t\ge 0$. Note that the positions $X_t$ and co-states $P_t$ may still undergo a time evolution; but need to satisfy
\begin{equation} \label{eq:steady state evolution}
    \dot{P}_t = D_x^2 \hat v (X_t) \dot{X}_t,
\end{equation}
which follows from (\ref{eq:PMP_forward}c) and $\nabla_x \phi_t = \nabla_x \hat v$. We are interested in the stationary solution $\nabla_x \hat v$ and run (\ref{eq:PMP_forward}) forward in time starting from appropriate initial conditions.
\end{remark}

Beyond steady state behavior, the following theorem characterizes the time evolution of the system (\ref{eq:PMP_forward}).

\begin{theorem}
    Assume that the equilibrium solution $\hat v$ of the stationary HJB equation
    (\ref{eq:stationary HJB}) is sufficiently smooth and that Assumption \ref{assumption4} holds.
    Given initial conditions (\ref{eq:ICmap}) for initial distribution $p_0>0$ and co-states $P_0 = 0$, the forward evolution equations (\ref{eq:PMP_forward}) lead to a time-evolving $\phi_t$ that satisfies (\ref{eq:HJB2d}) up to an irrelevant time-dependent constant; {\it i.e.}, $\nabla_x \phi_t = \nabla_x v_t$. Furthermore, $u_t = -RG^{\rm T}\nabla \phi_t$ and
    \begin{equation*}
    \lim_{t\to \infty} u_t = \hat u .
    \end{equation*}
\end{theorem}

\begin{proof}
Theorem \ref{theorem1} is also applicable to the specific choice (\ref{eq:time reversed beta}). The theorem therefore follows from Theorem \ref{theorem1} after time reversal. More specifically, we first note that $y_t = \nabla_x v_t$
satisfies
\begin{equation*}
\partial_t y_t = -\gamma y_t + \nabla_x \langle b,y_t \rangle - \frac{1}{2}\nabla_x \|RG^{\rm T}y_t\|^2_R +
    \frac{1}{2} \nabla_x (\Sigma : D_x y_t) + \nabla_x c , \qquad
    y_T = 0.
\end{equation*}
Next, we demonstrate that $\nabla_x \phi_t$ satisfies the same evolution equation:
\begin{align*}
    \partial_t \nabla_x \phi_t(X_t) =& \,-D_x^2 \phi_t(X_t)\dot{X}_t + \dot{P}_t\\ \nonumber
        =& \,-D_x^2 \phi_t(X_t) \dot{X}_t -\gamma P_t + \left(D_x (b(X_t)+G(X_t)U_t)\right)^{\rm T}P_t + \nabla_x c(X_t)
        \\ \nonumber
        &+ \frac{1}{2} \nabla_x  (\Sigma : D_x^2\phi_t)(X_t)  
        + D_x^2 \phi_t(X_t)\left(2\dot{X}_t+ \frac{1}{2}\Sigma \nabla_x \log \rho_t(X_t)\right)  \\ \nonumber
        =& \,-\gamma \nabla_x \phi_t(X_t) + \nabla_x \bigl< b(X_t),\nabla_x \phi_t(X_t)\bigr>  -\frac{1}{2} \nabla_x \|R G(X_t)^{\rm T}\nabla_x \phi_t(X_t)\|^2_R 
        \\
        & \,+ \frac{1}{2} \nabla_x  (\Sigma : D_x^2 \phi_t)(X_t) + \nabla_x c(X_t).
\end{align*}
Since the equality is valid for all $X_t$, the desired equivalence $\nabla_x \phi_t = y_t$ follows.
\end{proof}

We close this section by mentioning that (\ref{eq:PMP_forward}) is related to iterative approaches for approximating the stationary value function (or $Q$-function) used in reinforcement learning and based on the stationary Bellman equation \cite{Meyn}. Our approach instead targets the gradients of the value function (momenta/co-states) directly. This has been found to be beneficial from a theoretical perspective \cite{Carmona} and provides flexibility in terms of constraints added to the states and controls. It is also often found that approximating the gradient directly instead of its underlying  function is more robust computationally. Compare, for example, density versus score estimation \cite{MRO20,GLMRT25}.

%
\section{Linear control problems} \label{sec:linear}
%

We analyze the proposed formulations in the context of controlled linear SDEs
\begin{equation*}
{\rm d}X_t = AX_t + G U_t + \Sigma^{1/2}{\rm d}B_t, \qquad X_0 \sim {\rm N}(0,C_0)
\end{equation*}
for infinite time discounted cost functional (\ref{eq:costJ infinite}) with running cost
\begin{equation*}
c(x) = \frac{1}{2}x^{\rm T}C x
\end{equation*}
and $R=I$. We write the function $\phi_t(x)$ as
\begin{equation*}
    \phi_t(x) = \frac{1}{2}x^{\rm T}\Omega_t x.
\end{equation*}
With this {\it ansatz} the evolution equations (\ref{eq:time-reversed PMP}a)-(\ref{eq:time-reversed PMP}b) take the form
\begin{subequations}
    \begin{align*}
        \dot{X}_t &= AX_t - G G^{\rm T}P_t - \beta_t,\\
        -\dot{P}_t &= -\gamma P_t + A^{\rm T}P_t + CX_t + \Omega_t \beta_t.
    \end{align*}
\end{subequations}
The choice (\ref{eq:beta natural}) leads to
\begin{equation*}
    \beta_t = -\frac{1}{2}\Sigma ({\rm C}_t^{xx})^{-1}X_t
\end{equation*}
with ${\rm C}_t^{xx}$ the covariance matrix of the Gaussian distributed $X_t$ with mean zero. The constraint (\ref{eq:time-reversed PMP}c) becomes
\begin{equation*}
P_t = \Omega_t X_t,
\end{equation*}
which yields
\begin{subequations}
    \begin{align*}
        -\dot{\Omega}_t X_t &=  \Omega_t \dot{X}_t-\dot{P}_t\\
        &=\Omega_t \left(AX_t - G G^{\rm T}P_t - \beta_t\right)
        -\gamma P_t + A^{\rm T}P_t + CX_t + \Omega_t \beta_t\\
        &= \left(-\gamma \Omega_t + \Omega_t A - \Omega_t G G^{\rm T}\Omega_t + A^{\rm T}\Omega_t + C \right) X_t
    \end{align*}
\end{subequations}
and from which we can infer the backward Riccati equation for $\Omega_t$ \cite{Carmona}:
\begin{equation*}
-\dot{\Omega} = -\gamma \Omega_t + \Omega_t A + A^{\rm T}\Omega_t - \Omega_t G G^{\rm T}\Omega_t  + C.
\end{equation*}

Let us now investigate the forward evolution equations (\ref{eq:PMP_forward}). We obtain
\begin{subequations}
    \begin{align*}
        \dot{X}_t &= AX_t - G G^{\rm T}P_t - \widetilde{\upbeta}_t\\
        \dot{P}_t &= -\gamma P_t + A^{\rm T}P_t + CX_t +  \Omega_t \widetilde{\upbeta}_t 
        + 2 \Omega_t \dot{X}_t
    \end{align*}
\end{subequations}
subject to the constraint $P_t = \Omega_t X_t$. This time we obtain
\begin{subequations}
    \begin{align*}
        \dot{\Omega}_t X_t &=  -\Omega_t \dot{X}_t+\dot{P}_t\\
        &= \Omega_t \left(AX_t - G G^{\rm T} P_t - \beta_t\right)
        - \gamma P_t + A^{\rm T} P_t + CX_t +  \Omega_t \beta_t \\
        &= \left(- \gamma \Omega_t + \Omega_t A  + A^{\rm T}\Omega_t - \Omega_t G G^{\rm T}\Omega_t + C \right) X_t
    \end{align*}
\end{subequations}
and infer the forward Riccati equation
\begin{equation*}
\dot{\Omega} = -\gamma \Omega_t + \Omega_t A + A^{\rm T}\Omega_t - \Omega_t G G^{\rm T}\Omega_t  + C.
\end{equation*}
Clearly, the forward Riccati equation is stable if and only if the backward Riccati equation is stable and their stationary solutions agree. Of course, our approach is to propagate $X_t$ and $P_t$ and then find the appropriate $\Omega_t$ by enforcing the constraint $P_t = \Omega_t X_t$ through regression; {\it i.e.}
\begin{equation*}
\Omega_t = \arg \min_{\Omega\in \mathbb{R}^{d_x\times d_x}} \mathbb{E}\left[ \|P_t - \Omega X_t\|^2\right].
\end{equation*}
We describe numerical implementations for the general nonlinear setting in the following section.

%
\section{Numerical implementation} \label{sec:implementation}
%

Numerical implementations of the McKean--Pontryagin formulations rely on a set of particles $(X_t^{(i)},P_t^{(i)})$, $i=1,\ldots,M$, with corresponding evolution equations, which approximate the infinite-dimensional Hamiltonian system (\ref{eq:Hamiltonian_Lagrange}). At this point, two principal choices arise; either one discretizes directly at the level of (\ref{eq:Hamiltonian_Lagrange}) or, alternatively, one discretizes at the level of the Hamiltonian (\ref{eq:Hamiltonian_function_Lagrange}) first and then derives the corresponding finite-dimensional constrained Hamiltonian equations of motion. Both approaches will be discussed below. More specifically, we discuss three different particle-based approximations, which mainly differ in the manner in which the constraint (\ref{eq:Hamiltonian_Lagrange}c) is dealt with. These approximations apply equally well to the discounted cost formulation (\ref{eq:PMP_forward}) of the McKean--Pontryagin minimum principle from Section \ref{sec:discounted}.

\subsection{A regression formulation}
The constraint (\ref{eq:Hamiltonian_Lagrange}c) can be approximated using the least squares method. For example, one introduces the cost function
\begin{equation} \label{eq:least squares}
\mathcal{V}_t(\eta) = \frac{1}{2M} \sum_{i=1}^M
\| P_t^{(i)} - \tilde \Phi(X_t^{(i)};\eta)\|^2
\end{equation}
for a suitable class of functions $\tilde \Phi(\cdot\,;\eta)$ parametrized by $\eta \in \mathbb{R}^{d_\eta}$. The approximation $\tilde \Phi(\cdot\,;\eta_t)$ with $\eta_t = 
\arg \min \mathcal{V}_t(\eta)$ is then used in
(\ref{eq:Hamiltonian ODE}) instead of $\nabla_x \phi_t(x)$ and gives rise to associated equations of motion in terms of the particles $(X_t^{(i)},P_t^{(i)})$, $i=1,\ldots,M$. 
 
We note that a regression approach has been combined with the classical PMP for controlled ordinary differential equations in \cite{Kunisch23} by including a penalty of the form (\ref{eq:least squares}) in the cost function (\ref{eq:costJ}).

Although the regression approach can be combined with any suitable choice of $\upbeta_t$, the following two approaches are based on the specific choice (\ref{eq:beta natural}).

%
\subsection{Discrete Schr\"odinger bridge approximation} 
\label{sec:DM}
%

Consider in more detail the formulation (\ref{eq:Hamiltonian ODE}) with $\widetilde{\upbeta}_t$ given by (\ref{eq:beta natural}). We replace the pull-back $\widehat{\mathcal{L}}_{\widetilde{\beta}_t}$ of the generator (\ref{eq:natural generator}) by a discrete Schr\"odinger bridge approximation $m_t^{(ij)}$ \cite{MarshallCoifman,WR20,GLRY24} over the set of particles $\{X_t^{(j)}\}$  giving rise to the following particle approximation of the Hamiltonian (\ref{eq:Hamiltonian_function_Lagrange_generator_simple}):
\begin{subequations}
\begin{align*}
\widetilde{\mathrm{H}}_M(\{X_t^{(j)}\},\{P_t^{(j)}\},\{U_t^{(j)}\}) =&\,\, \frac{1}{M}
\sum_{i=1}^M
\left\{ \bigl< b(X_t^{(i)}) + G(X_t^{(i)})U_t^{(i)},P_t^{(i)} \bigr> 
+ \frac{1}{2} \|U_t^{(i)}\|^2\right\} \\
&+ \frac{1}{M}\sum_{i=1}^M \left\{
c(X_t^{(i)}) - \sum_{j=1}^M m_t^{(ij)} \bigl< X_t^{(j)},P_t^{(i)}\bigr>  \right\}.
 \end{align*}
\end{subequations}
The coefficients $m_t^{(ij)}$ are defined by
\begin{equation*}
    m_t^{(ij)} = \frac{1}{\varepsilon} \left( v_t^{(i)} d_t^{(ij)} v_t^{(j)} - \updelta_{ij}\right), \qquad d_t^{(ij)} = 
    e^{-\frac{1}{2\varepsilon}\|X_t^{(i)}-X_t^{(j)}\|_\Sigma^2},
\end{equation*}
with $\{v_t^{(i)}\ge 0\}$ chosen such that 
\begin{equation*}
\sum_i m_t^{(ij)} = \sum_j m_t^{(ij)} = 0.
\end{equation*}
Here $\updelta_{ij}$ denotes the Kronecker delta function; {\it i.e.}, $\updelta_{ij} = 1$ if $i=j$ and $\updelta_{ij} = 0$ otherwise. The regularization parameter $\varepsilon>0$ should be chosen appropriately \cite{WR20}. In \cite{MarshallCoifman,WR20} a fast iterative algorithm for computing coefficients $v_t^{(i)}$ has been introduced.

The associated Hamiltonian equations of motion are given by
\begin{subequations} \label{eq:interacting Hamiltonian ODE}
    \begin{align}
        M^{-1}\dot{X}_t^{(i)} &= \nabla_{P^{(i)}} \widetilde{\mathrm{H}}_M(\{X_t^{(j)}\},\{P_t^{(j)}\},\{U_t^{(j)}\}),\\
        -M^{-1} \dot{P}_t^{(i)} &= \nabla_{X^{(i)}} \widetilde{\mathrm{H}}_M(\{X_t^{(j)}\},\{P_t^{(j)}\},\{U_t^{(j)}\})\\
        0&= G(X_t^{(i)})^{\rm T}P_t^{(i)} + R^{-1}U_t^{(i)}.
    \end{align}
\end{subequations}
$i=1,\ldots,M$, subject to boundary conditions
\begin{equation*}
X_0^{(i)}  \sim p_0, \qquad
P_T^{(i)} = \nabla_x f(X_T^{(i)}).
\end{equation*}
The additional $M^{-1}$ term arises from replacing Frechet derivatives $\delta \cdot/\delta X_t$ and $\delta \cdot /\delta P_t$ by gradients with respect particle positions $X_t^{(i)}$ and co-states $P_t^{(i)}$, respectively. 

Computing gradients of $m_t^{(ij)}$ with respect to particle positions $X_t^{(k)}$ is non-trivial, and we use the simpler approximation based on the following direct application of the Schr\"odinger bridge approximation to (\ref{eq:Hamiltonian ODE}):
\begin{subequations}
    \begin{align*}
        \dot{X}_t^{(i)} =& \,\,b(X_t^{(i)}) - G(X_t^{(i)})R G(X_t^{(i)})^{\rm T} P_t^{(i)}
        - \sum_{j=1}^M m_t^{(ij)}X_t^{(j)},\\
        -\dot{P}_t^{(i)} =& \,\,(D_xb(X_t^{(i)}))^{\rm T}P_t^{(i)} +\nabla_x c(X_t^{(i)})\,
         \\
        &  \,-\,\frac{1}{2} \nabla_x 
        \|RG(X_t^{(i)})^{\rm T}P_t^{(i)}\|^2_R + \sum_{j=1}^M m_t^{(ij)}P_t^{(j)}.
    \end{align*}
\end{subequations}
The discrete Schr\"odinger bridge approximation does not require an explicit approximation of $\phi_t(x)$, which simplifies numerical implementations. Time-stepping with step-size $\Delta t>0$ is achieved by a combination of forward and backward Euler discretization:
\begin{subequations}\label{eq:time-stepping FB Euler}
    \begin{align} 
        \frac{X_{t_{n+1}}^{(i)} -X_{t_n}^{(i)} }{\Delta t}=& \,\,b(X_{t_n}^{(i)}) - G(X_{t_n}^{(i)})R G(X_{t_n}^{(i)})^{\rm T} P_{t_n}^{(i)}
        - \sum_{j=1}^M m_{t_n}^{(ij)}X_{t_n}^{(j)},\\ \nonumber
        \frac{P_{t_{n}}^{(i)} -P_{t_{n+1}}^{(i)} }{\Delta t}
         =& \,\,(D_xb(X_{t_{n+1}}^{(i)}))^{\rm T}P_{t_{n+1}}^{(i)} +\nabla_x c(X_{t_{n+1}}^{(i)})\,
         \\
        &  \,-\,\frac{1}{2} \nabla_x 
        \|RG(X_{t_{n+1}}^{(i)})^{\rm T}P_{t_{n+1}}^{(i)}\|^2_R + \sum_{j=1}^M m_{t_{n+1}}^{(ij)}P_{t_{n+1}}^{(j)},
    \end{align}
\end{subequations}
$t_{n+1} = t_n+\Delta t$, $n=0,\ldots,N-1$, $T = \Delta t\,N$.

The corresponding interacting particle equations for infinite horizon discounted optimal control problems are provided by
\begin{subequations} \label{eq:forward formulation}
    \begin{align}
        \dot{X}_t^{(i)} =& \,\,b(X_t^{(i)}) - G(X_t^{(i)})R G(X_t^{(i)})^{\rm T} P_t^{(i)}
        - \sum_{j=1}^M m_t^{(ij)}X_t^{(j)},\\ \nonumber
        \dot{P}_t^{(i)} =& \,\,-\gamma P_t + (D_xb(X_t^{(i)}))^{\rm T}P_t^{(i)} +\nabla_x c(X_t^{(i)})\,
         \\
        &  \,-\,\frac{1}{2} \nabla_x 
        \|RG(X_t^{(i)})^{\rm T}P_t^{(i)}\|^2_R + \sum_{j=1}^M m_t^{(ij)}P_t^{(j)} +
        2 D_x^2 \phi_t(X_t^{(i)})\dot{X}_t^{(i)}.
    \end{align}
\end{subequations}
Equations (\ref{eq:forward formulation}) are time-stepped using forward Euler for both $X_t^{(i)}$ and $P_t^{(i)}$. We use this approach for numerical experiments in Sections \ref{sec:num_IP1} and \ref{sec:num_L63}.

%
\subsection{Variational approach} \label{sec:variational}
%

A particularly versatile approach is to approximate $\upbeta_t$ and $\phi_t$ via a parametric {\it ansatz} $\upbeta_t(x) = \tilde \upbeta(x;\eta_t)$ and $\nabla_x \phi_t(x) = \tilde \Phi(x;\alpha_t)$, respectively, and to consider the Hamiltonian
\begin{align} \nonumber
    &\mathrm{H}_M(\{X_t^{(j)}\},\{P_t^{(j)}\},\{U_t^{(j)}\},\eta_t,\alpha_t) \,=\\  \nonumber
    &\qquad \qquad \frac{1}{M} \sum_{i=1}^M \left\{
    \bigl< b(X_t^{(i)}) + G(X_t^{(i)})U_t^{(i)},P_t^{(i)} \bigr> 
+ \frac{1}{2} \|U_t^{(i)}\|^2   + c(X_t^{(i)})  \right. \\
& \qquad \qquad \left. + \tilde \upbeta(X_t^{(i)},\eta_t) \cdot  \left(\tilde \Phi(X_t^{(i)},\alpha_t)- P_t^{(i)}\right) +
\frac{1}{2} \nabla_x \cdot \left( \Sigma \tilde \Phi(X_t^{(i)};\alpha_t) \right) \right\}
\label{eq:particle Hamiltonian}
\end{align}
with equations of motion
\begin{subequations} \label{eq:variational ODE}
    \begin{align}
        M^{-1} \dot{X}_t^{(i)} &= 
        \nabla_{P^{(i)}}\mathrm{H}_M(\{X_t^{(j)}\},\{P_t^{(j)}\},\{U_t^{(j)}\},\eta_t,\alpha_t),\\
        -M^{-1} \dot{P}_t^{(i)} &= \nabla_{X^{(i)}}\mathrm{H}_M(\{X_t^{(j)}\},\{P_t^{(j)}\},\{U_t^{(j)}\},\eta_t,\alpha_t),\\
        0 &= \nabla_\eta \mathrm{H}_M(\{X_t^{(j)}\},\{P_t^{(j)}\},\{U_t^{(j)}\},\eta_t,\alpha_t)\\
        0 &= \nabla_\alpha \mathrm{H}_M(\{X_t^{(j)}\},\{P_t^{(j)}\},\{U_t^{(j)}\},\eta_t,\alpha_t), \\
        0&= G(X_t^{(i)})^{\rm T}P_t^{(i)} + R^{-1}U_t^{(i)}.
    \end{align}
\end{subequations}
Time-stepping is achieved with the same combination of forward and backward Euler approximations as used in (\ref{eq:time-stepping FB Euler}). Suitable choices for $\tilde \upbeta$ and $\tilde \Phi$ can follow the approaches considered, for example, in \cite{MRO20}. In the following, we discuss the case of linear approximations in some detail.

\begin{example}
    Consider the linear {\it ansatz}
    \begin{equation*}
    \tilde \upbeta(x;\eta) = B x + d, \qquad
    \tilde \Phi(x;\alpha) = Ax + c
    \end{equation*}
    for suitable matrices $A, B$ and vector $c, d$; {\it i.e.}, $\eta = (d,B)$ and $\alpha = (c,A)$. Taking variations with respect to these parameters in (\ref{eq:particle Hamiltonian}) leads to 
    \begin{equation} \label{eq:linear regression}
    A_t = {\rm C}_t^{px}({\rm C}_t^{xx})^{-1}, \qquad
    c_t = \mu^p_t-A_t \mu^x_t
    \end{equation}
    and
    \begin{equation*}
        B_t = -\frac{1}{2} \Sigma ({\rm C}_t^{xx})^{-1}, \qquad d_t = -B_t \mu_t^x,
    \end{equation*}
    respectively. Here, ${\rm C}_t^{px}$ and ${\rm C}_t^{xx}$ denote the empirical covariance matrices between $\{P_t^{(i)}\}$ and $\{X_t^{(i)}\}$, and $\{X_t^{(i)}\}$ with itself, respectively. The empirical mean of $\{X_t^{(i)}\}$ and $\{P_t^{(i)}\}$, respectively, is denoted by $\mu_t^x$ and $\mu_t^p$, respectively. Upon substituting these parameters into (\ref{eq:particle Hamiltonian}), the following Hamiltonian equations of motion arise from (\ref{eq:interacting Hamiltonian ODE}):
    \begin{subequations} \label{eq:linear variational}
    \begin{align}
        \dot{X}_t^{(i)} = & \,\,b(X_t^{(i)}) - G(X_t^{(i)})R G(X_t^{(i)})^{\rm T} P_t^{(i)} + \frac{1}{2} \Sigma ({\rm C}_t^{xx})^{-1}(X_t^{(i)}-\mu_t^x),\\ \nonumber
        -\dot{P}_t^{(i)} = & \,\,(D_x b(X_t^{(i)}))^{\rm T}P_t^{(i)} +\nabla_x c(X_t^{(i)})
        - \frac{1}{2} \nabla_x 
        \|RG(X_t^{(i)})^{\rm T}P_t^{(i)}\|^2_R  \\\nonumber & 
        + \frac{1}{2} ({\rm C}_t^{xx})^{-1}\Sigma \left( P_t^{(i)}-\mu_t^p - {\rm C}^{px}_t ({\rm C}_t^{xx})^{-1}(X_t^{(i)}-\mu_t^x) 
        \right) \\ &
        - \frac{1}{2} ({\rm C}_t^{xx})^{-1}{\rm C}_t^{xp}  \Sigma ({\rm C}_t^{xx})^{-1}(X_t^{(i)}-\mu_t^x).
    \end{align}
    \end{subequations}
    Note that derivatives of the parameters with respect to the particle positions and momenta do not appear because of (\ref{eq:variational ODE}c)--(\ref{eq:variational ODE}d). Such a linear approximation of the McKean--Pontryagin minimum principle is of interest for receding horizon stochastic optimal control problems as they arise from NMPC. See Section \ref{sec:num_IP2} for a numerical illustration.
\end{example}

We close this section by mentioning that the deep learning approaches exploited in 
\cite{EHJ17,EHJ21,nusken2020solving,hua2025an} can also be extended to our McKean--Pontryagin formulations. In this case, one would utilize global parametric {\it ansatz} functions $\upbeta_t(x) = \tilde \upbeta_t(x;\eta)$ and $\nabla_x \phi_t(x) = \tilde \Phi_t(x;\alpha)$, $t\in [0,T]$, and directly substitute them into the action principle (\ref{eq:action Lagrange}).

%
\section{Numerical experiments} \label{sec:experiments}
%

We present results from five numerical experiments. The first experiment considers a controlled one-dimensional diffusion process for which a reference solution can be computed by numerical approximation of the corresponding HJB equation. The numerical convergence of an interacting particle approximation of (\ref{eq:NMPC formulation}) is demonstrated. The remaining four experiments are qualitative and demonstrate the broad scope of the proposed methodology. The second and third experiment involve infinite horizon cost functions of the form (\ref{eq:costJ infinite}) and controlled ODEs. In both cases, we add a diffusion (viscosity) regularization in (\ref{eq:Hamiltonian_Lagrange}) in line with the concept of viscosity solutions to the associated HJB equations and the vanishing viscosity method \cite{CL83}. The fourth experiment returns to the inverted pendulum from the second experiment but now in the context of NMPC and receding horizon optimal control. Here, we rely on an implementation of the formulation (\ref{eq:linear variational}). Our final example demonstrates the applicability of our approach to a higher dimensional $(d_x=40)$ Lorenz-96 system \cite{lorenz96}. The goal is to stabilize the (unstable) reference solution $x_{\rm ref} = 2$. We formulate the task as an infinite horizon control problem and use the linear approximation (\ref{eq:linear regression}) for $\nabla_x \phi_t(x)$.

%
\subsection{Diffusion under a bistable potential: Finite horizon control}
%

\begin{figure}
\centerline{\includegraphics[width=0.5\textwidth]{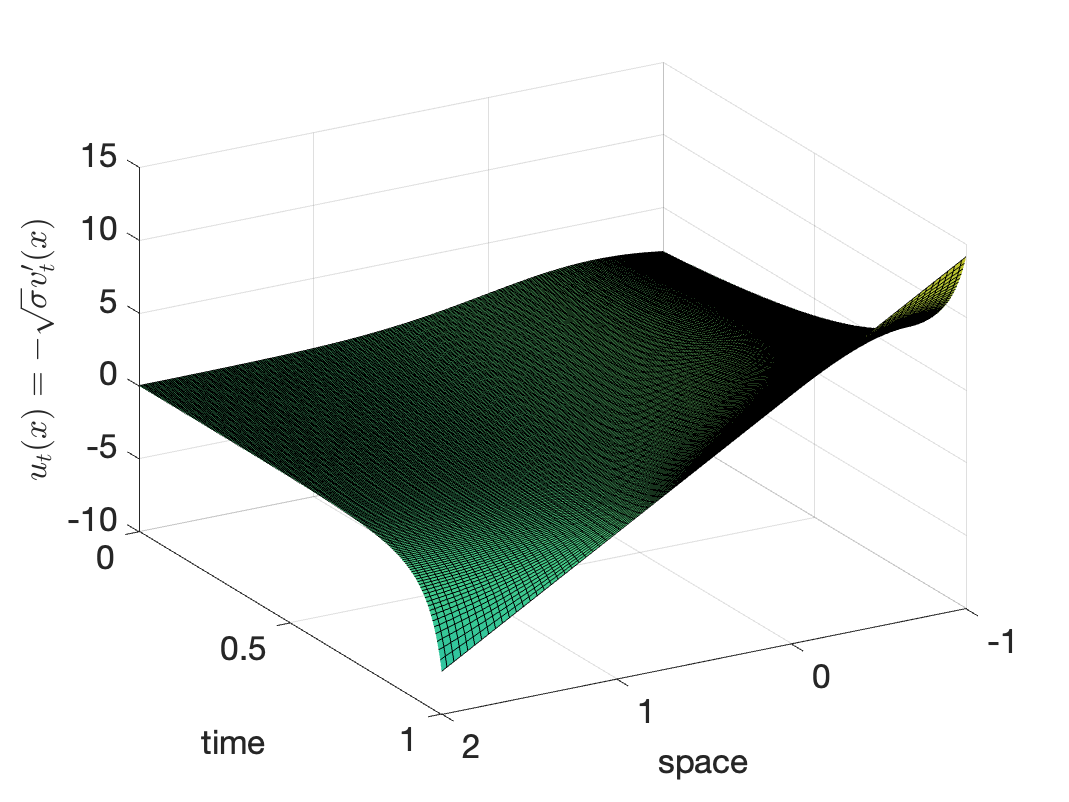}
$\quad$\includegraphics[width=0.42\textwidth]{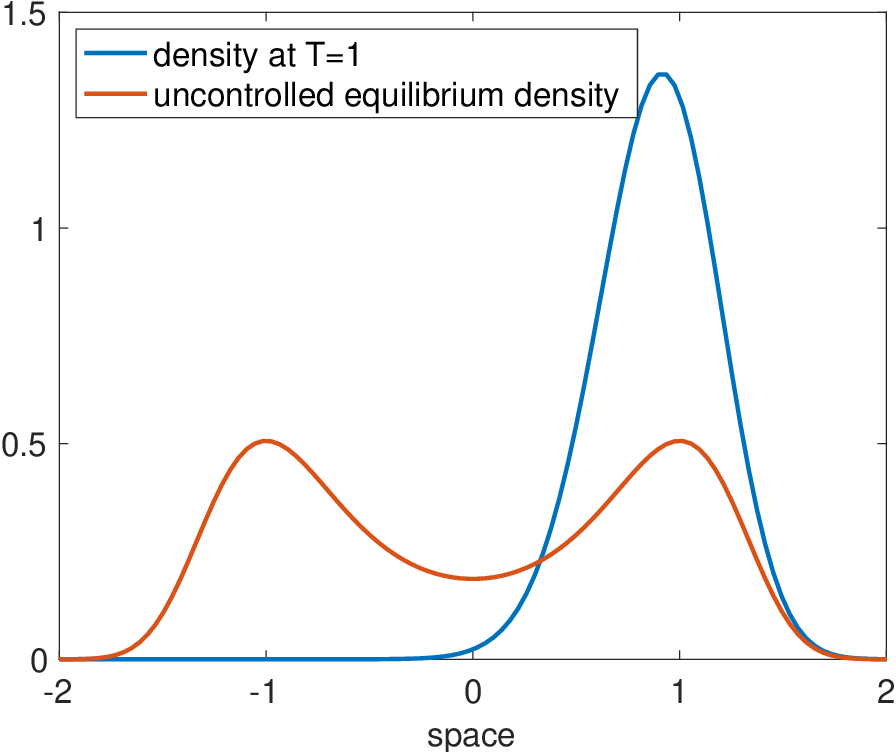}}
\caption{Left panel: Control $u_t(x)$ as a function of space and time. Right panel: Distribution of $X_T$ of the controlled process at final time $T=1$ and equilibrium distribution of the uncontrolled process.
} 
\label{fig:ex1a}
\end{figure}

We consider a one-dimensional controlled diffusion process
\begin{equation} \label{eq:FSDE ex1}
    {\rm d}X_t = -V^\prime(X_t)\,{\rm d}t + \sqrt{\sigma}\left(
    U_t\,{\rm d}t + {\rm d}B_t\right)
\end{equation}
with double-well potential
\begin{equation*}
    V(x) = \frac{1}{4} x^4 - \frac{1}{2}x^2
\end{equation*}
and diffusion constant $\sigma = 1/2$. The control task is to stabilize
$x_{\rm eq} = 1$ by minimizing the finite horizon cost (\ref{eq:costJ}) for $x_0 =0$, $T=1$, $R=1$, $c(x) = 0$, and terminal cost
\begin{equation*}
    f(x) = 5(x-1)^2.
\end{equation*}

We first compute the feedback control law $u_t(x)$ by solving the (linear) 
HJB equation in $w_t(x) = \exp(-v_t(x))$ numerically in a spatial domain $x\in [-4,4]$ with mesh-size $\Delta x = 0.04$ and step-size $\Delta t = 0.0002$. The resulting feedback control law
\begin{equation} \label{eq:OC_ex1}
u_t(x) = -\sqrt{\sigma}v'_t(x),
\end{equation}
the density $\rho_T(x)$ of $X_T$ for the controlled process with
$U_t=u_t(X_t)$, and the equilibrium density $\rho_{\rm eq}(x) \propto \exp(-2/\sigma V(x))$ of the uncontrolled process can be found in Figure \ref{fig:ex1a}. The value of
the associated finite horizon cost (\ref{eq:costJ}) with control (\ref{eq:OC_ex1}) is $v_0(0) 
\approx 1.4319$.

\begin{figure}
\centerline{\includegraphics[width=0.8\textwidth]{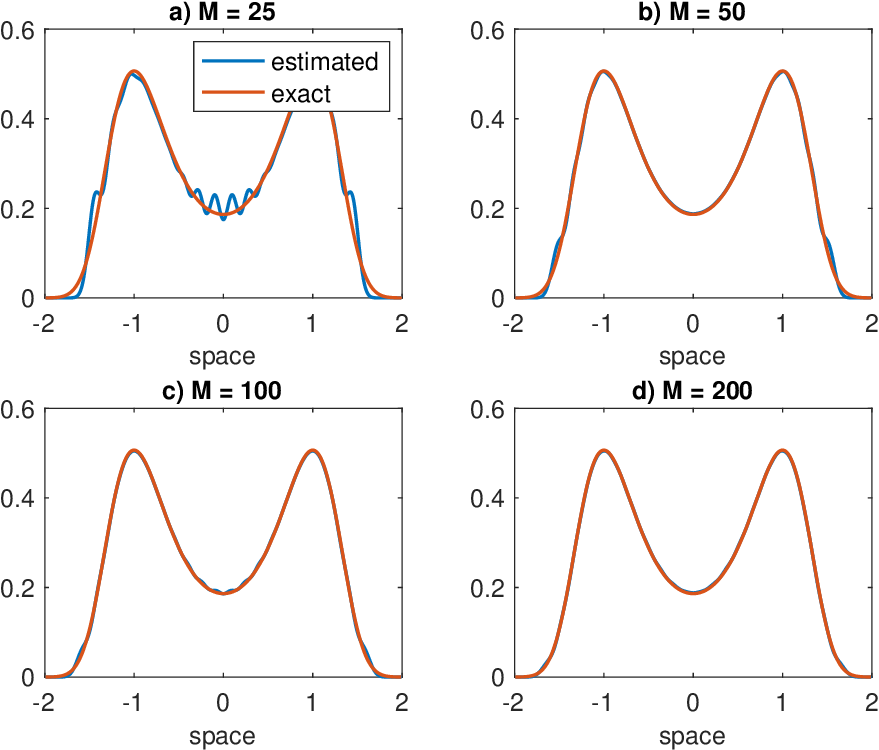}}
\caption{Displayed are the exact equilibrium and the 
approximated  distributions at time $t=4$ for fixed $\epsilon = 0.02$ and varying ensemble sizes: (a) $M=25$, (b) $M=50$, (c) $M=100$, and 
d) $M=200$.} 
\label{fig:ex1b}
\end{figure}

Next we approximate the invariant distribution $\rho_{\rm eq}(x)$ by solving the
uncontrolled mean-field ODE formulation
\begin{equation*}
    \dot{X}_t = -V'(X_t) + \widehat{\mathcal{L}}_{\widetilde{\beta_t}}X_t
\end{equation*}
numerically over the time interval $t\in [0,4]$ for ensemble sizes $M = 25,50,100$, and $200$, step-size $\Delta t = 0.01$, and initial condition $X_0 \sim {\rm N}(0,0.01)$ The operator $\widehat{\mathcal{L}}_{\widetilde{\beta_t}}$ is approximated numerically using discrete Schr\"odinger bridges as outlined in Section \ref{sec:DM}.  See Table \ref{table1} for approximation errors for varying
parameter values $\epsilon$ and ensemble sizes $M$. Figure \ref{fig:ex1b} displays density approximations for $\epsilon = 0.02$ and varying ensemble sizes $M$. We employ $\epsilon = 0.02$ also in the subsequent simulations.

\begin{table}
\begin{center}
\begin{tabular}{ |c|c|c|c| } 
 \hline
$\epsilon/M$  & $0.01$ & $0.02$ & $0.03$  \\ 
\hline 
 $25$ & 0.0507  & 0.0257 & 0.0311\\ 
 $50$ & 0.0109  & 0.0078 & 0.0099\\
 $100$ & 0.0072  & 0.0044 & 0.0076\\
 $200$ & 0.0059  & 0.0028 & 0.037\\
 \hline
\end{tabular}
\end{center}
\caption{Approximation errors in TV norm for the computed approximations of the equilibrium density $\rho_{\rm eq}$ as a function of the ensemble size $M$ and the parameter $\epsilon$.}
\label{table1}
\end{table}

Finally, we employ formulation (\ref{eq:NMPC formulation}) for four different configurations; {\it i.e.}, with reference control $u_t^{\rm ref}(x) \equiv 0$ and 
\begin{equation} \label{eq:linear ref}
u_t^{\rm ref}(x) = -(x-1),
\end{equation}
respectively, and initial distributions 
\begin{equation}
    p_0(x) = \frac{1}{\sqrt{2\pi s}}e^{\frac{-1}{2s}x^2}
\end{equation}
with $s=1$ and $s = 0.01$, respectively. The particle and temporal discretizations follow
(\ref{eq:time-stepping FB Euler}). The additional advection term in (\ref{eq:NMPC formulation}b) is approximated using Nadaraya--Watson kernel regression \cite{Bierens} for the constraint (\ref{eq:NMPC formulation}c); {\it i.e.},
\begin{equation} \label{eq:NWKR}
\nabla_x \phi_t(x) \approx \breve{y}_t(x):=
\frac{\sum_{i=1}^M \exp \left(-\frac{1}{2\delta}\|x - X_t^{(i)}\|^2 \right) P_t^{(i)}}{\sum_{i=1}^M \exp \left(-\frac{1}{2\delta}\|x - X_t^{(i)}\|^2 \right)}
\end{equation}
for $\delta = \epsilon$, and
\begin{equation} \label{eq:NWKR_FD}
    D_x^2 \phi_t (X_t)\,{\rm b}_t(X_t,P_t) \approx \frac{1}{\Delta t} \left(\breve{y}_t(X_t+\Delta t \,{\rm b}_t(X_t,P_t)) - \breve{y}_t(X_t)\right)
\end{equation}
with
\begin{equation*}
    {\rm b}_t(X_t,P_t) = G(X_t)\left(RG(X_t)^{\rm T} P_t +  u_t^{\rm ref}(X_t)\right).
\end{equation*}

\begin{figure}
\centerline{\includegraphics[width=0.8\textwidth]{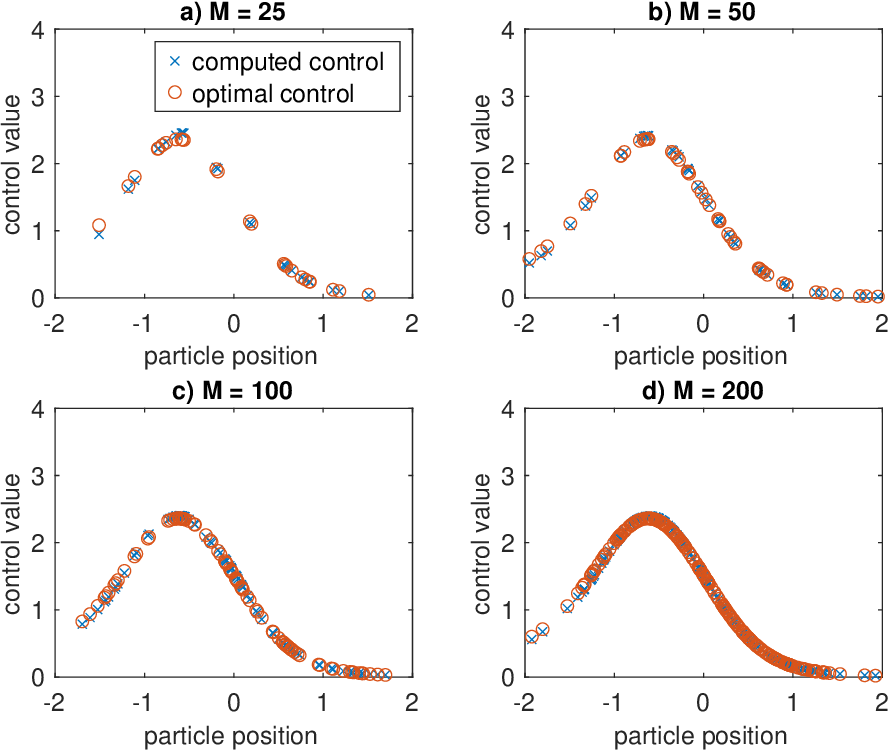}}
\caption{Displayed are the optimal control $u_0(X_0^{(i)})$, as computed from the HJB equation, at particle positions $X_0^{(i)}$ and the control $U_0^{(i)} = -
\sqrt{\sigma}P_0^{(i)}$ computed by the McKean--Pontryagin formulation at time $t=0$ for varying ensemble sizes $M$. The reference control is set to zero and the variance of the initial particle distribution is $s=1$.} 
\label{fig:ex1c}
\end{figure}

\begin{figure}
\centerline{\includegraphics[width=0.8\textwidth]{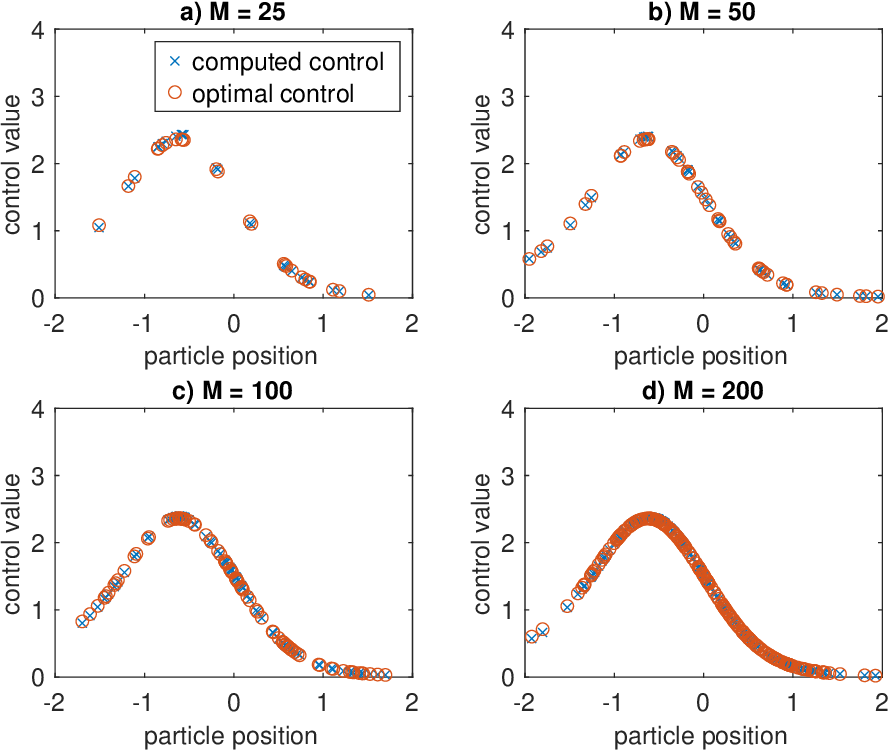}}
\caption{Displayed are the optimal control $u_0(X_0^{(i)})$, as computed from the HJB equation, at particle positions $X_0^{(i)}$ and the control $U_0^{(i)} = -
\sqrt{\sigma}P_0^{(i)}$ computed by the McKean--Pontryagin formulation at time $t=0$ for varying ensemble sizes $M$. The reference control is set to $u_t^{\rm ref}(x) = -(x-1)$ and the variance of the initial particle distribution is $s=1$.} 
\label{fig:ex1d}
\end{figure}

\begin{figure}
\centerline{\includegraphics[width=0.8\textwidth]{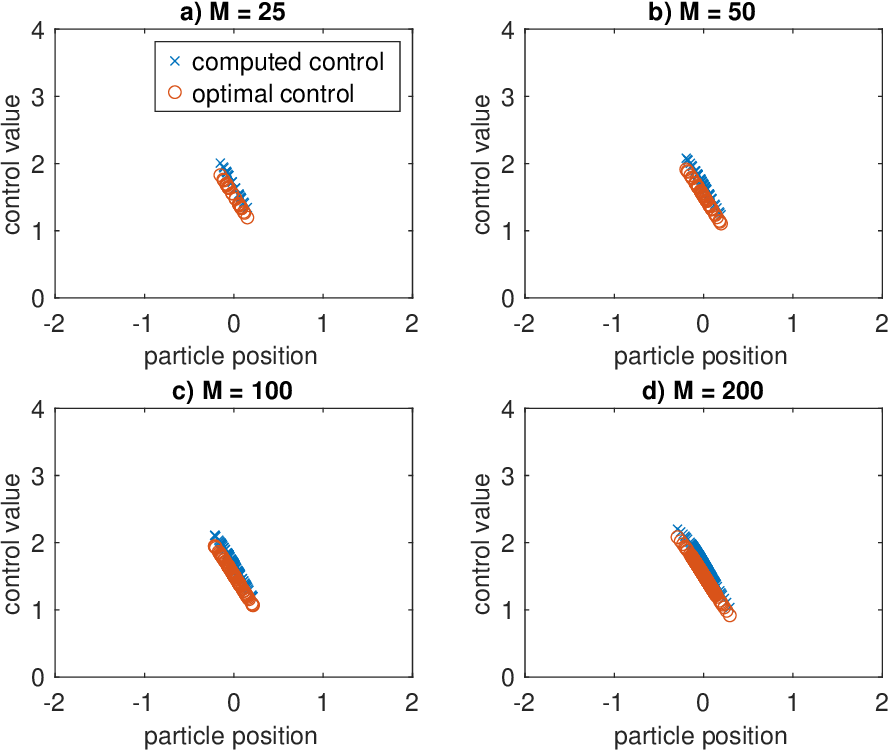}}
\caption{Displayed are the optimal control $u_0(X_0^{(i)})$, as computed from the HJB equation, at particle positions $X_0^{(i)}$ and the control $U_0^{(i)} = -
\sqrt{\sigma}P_0^{(i)}$ computed by the McKean--Pontryagin formulation at time $t=0$ for varying ensemble sizes $M$. The reference control is set to zero and the variance of the initial particle distribution is $s=0.01$.} 
\label{fig:ex1e}
\end{figure}

\begin{figure}
\centerline{\includegraphics[width=0.8\textwidth]{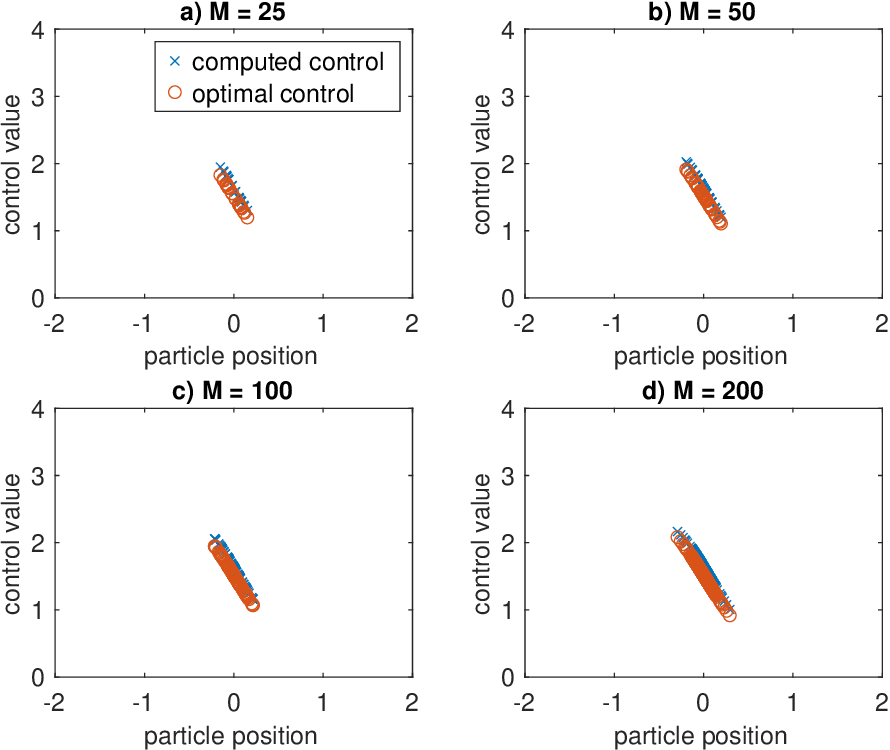}}
\caption{Displayed are the optimal control $u_0(X_0^{(i)})$, as computed from the HJB equation, at particle positions $X_0^{(i)}$ and the control $U_0^{(i)} = -
\sqrt{\sigma}P_0^{(i)}$ computed by the McKean--Pontryagin formulation at time $t=0$ for varying ensemble sizes $M$. The reference control is set to $u_t^{\rm ref}(x) = -(x-1)$ and the variance of the initial particle distribution is $s=0.01$.} 
\label{fig:ex1f}
\end{figure}

In a first step, the particles $\{X_{t_n}^{(j)}\}$ are propagated forward in time using the control law $u_t^{\rm ref}(x)$, step-size $\Delta t = 1/N$, and initial conditions $X_0 \sim p_0$. In a second step, the co-states $\{P_{t_n}^{(j)}\}$ are propagated backward in time subject to terminal conditions $P_{t_N}^{(i)} = f'(X_{t_N}^{(i)})$, $i=1,\ldots,M$. The numerical results for variance $s = 0.01,\,1$ in the initial conditions, $u_t^{\rm ref}(x) = 0$, and linear reference control (\ref{eq:linear ref})
can be found in Figures \ref{fig:ex1c} to \ref{fig:ex1f}. It is found that even small a small ensemble size of $M=25$ reproduces the control computed from the HJB equation accurately along particle positions $\{X_0^{(i)}\}$. A lower value of $s$ appears advantageous while an impact of the reference control on the approximation quality is not detectable. 

We finally use the computed control laws to forward simulate the SDE (\ref{eq:FSDE ex1})
and approximate the cost (\ref{eq:costJ}) using $M = 10^6$ Monte Carlo samples. 
The results can be found in Table \ref{table2}. It can be concluded that the proposed formulation is insensitive to the specific choices of $s$ and $u_t^{\rm ref}$ and that the approximation quality is remarkably good across all selected ensemble sizes $M$. The fact that $M=25$ even almost always outperforms larger ensemble sizes must be attributed to the fact that parameters $\epsilon$, $\delta$, and $\Delta t$ have not been optimized relative to ensemble sizes $M$.

\begin{table}
\begin{center}
\begin{tabular}{ |c|c|c|c|c| } 
 \hline
  & $s=1$ & $s=1$ & $s=0.01$ &
  $s=0.01$\\ 
  & $u_t^{\rm ref}(x) =0$ & $u_t^{\rm ref}(x)=-(x-1)$ & $u_t^{\rm ref}(x) =0$& $u_t^{\rm ref}(x)=-(x-1)$ \\
\hline 
 $M=25$  & 1.4348  & 1.4340 & 1.4375 & 1.4374  \\ 
 $M=50$  & 1.4359  & 1.4361 & 1.4397 & 1.4394 \\
 $M=100$ & 1.4343  & 1.4341 & 1.4387 & 1.4380 \\
 $M=200$ & 1.4342  & 1.4341 & 1.4387 & 1.4380 \\
 \hline
\end{tabular}
\end{center}
\caption{Value of the finite horizon cost (\ref{eq:costJ}) for controls computed by four different configurations of the McKean--Pontryagin formulation. The reference value obtained for the control from the HJB equation is $1.4319$.}
\label{table2}
\end{table}

Let us briefly discuss the computational complexity of the proposed methodology in relation to a grid-based discretization of the underlying HJB equation. The proposed implementations of the interacting particle formulations scale quadratically in $M$ and linearly in the number of time-steps. Grid-based methods scale linearly in both the number of grid points and number of time-steps. However, explicit time-stepping methods require a step-size that is dictated by numerical stability, which in turn is determined by spatial resolution. Furthermore, while the required ensemble sizes can potentially become dimension independent for Monte Carlo like particle approximations as proposed here for stochastic optimal control, standard grid-based methods require an exponentially increasing number of grid points as the dimension of the HJB equation increases. See Section \ref{sec:num_L96} for a higher-dimensional test problem for which a direct disretization of the HJB equation is infeasible. A more detailed complexity study is left for future research. 

%
\subsection{Inverted pendulum: Infinite horizon control} \label{sec:num_IP1}
%

\begin{figure}
\centerline{\includegraphics[width=0.9\textwidth]{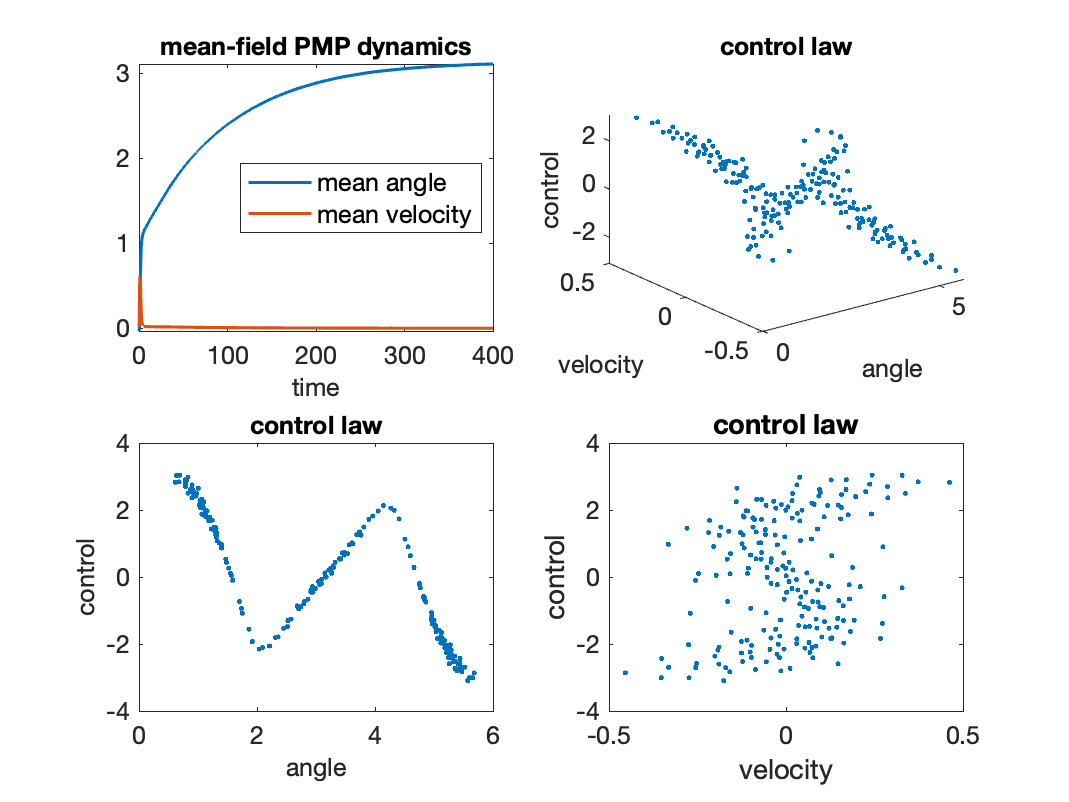}}
\caption{Upper left panel: Time evolution of the ensemble mean angle and ensemble mean velocity. Upper right panel: Control $U_T^{(i)}$
at final time $T=400$ as function of $X_T^{(i)} = 
(\theta_T^{(i)},v_T^{(i)})$, $i=1,\ldots,M$. Lower left panel: Control 
$U_T^{(i)}$ at final time as function of $\theta_T^{(i)}$.
Lower right panel: Control $U_T^{(i)}$ at final time as function of 
$v_T^{(i)}$.
} 
\label{fig1}
\end{figure}

We consider a controlled inverted pendulum with friction \cite{Meyn}. The state variable $x = (\theta,v)^{\rm T}$ is two-dimensional with equations of motion
\begin{subequations} \label{eq:pendulum}
\begin{align}
\dot{\theta}_t &= v_t,\\
\dot{v}_t &= -\sin (\theta_t) -\sigma v_t + \cos(\theta_t)U_t 
\end{align}
\end{subequations}
and $\sigma = 5$. Consider the running cost function 
\begin{equation} \label{eq:running cost IP}
c(x) = \frac{\alpha}{2} \left( (\theta-\pi)^2
+ v^2 \right)
\end{equation}
with $\alpha = 30$ and set $R=1$ in (\ref{eq:norm}). We use the discount factor $\gamma = 1.5$
in (\ref{eq:costJ infinite}). Note that $(\pi,0)^{\rm T}$ is an unstable equilibrium point under $U_t \equiv 0$.

In order to implement and, in fact, regularize the McKean--Pontryagin formulation, we add diffusion (viscosity) by setting $\Sigma = 0.1I$ in (\ref{eq:PMP_forward}). We initialize the particles at the stable equilibrium point $x_{\rm s} = (0,0)^{\rm T}$ with added Gaussian noise of variance $0.1I$; {\it i.e.}~$X_0^{(i)} \sim {\rm N}(x_{\rm s},0.1I)$. The initial momenta are  $P_0^{(i)} = (0,0)^{\rm T}$. The parameter $\delta$ in (\ref{eq:NWKR}) is set to $\delta = 0.1$ and 
${\rm b}_t = \dot{X}_t$ in (\ref{eq:NWKR_FD}).

The discrete Schr\"odinger bridge based formulation from Section \ref{sec:DM} has been implemented with ensemble size $M = 200$ and $\varepsilon = 2\Delta t$. The resulting evolution equations (\ref{eq:PMP_forward}) and its interacting particle approximation (\ref{eq:forward formulation}) have been simulated over 8000 forward Euler time-steps with time-step $\Delta t = 0.05$.

The results can be found in Figure \ref{fig1}. It can be clearly seen that the mean-field formulation leaves the stable equilibrium and rapidly equilibrates at the unstable equilibrium. The resulting final control terms 
\begin{equation} \label{eq:control law pendulum}
U_t^{(i)} = -G(\theta_t^{(i)})^{\rm T} P_t^{(i)}, \qquad G(\theta) = \left( \begin{array}{c} 0 \\ \cos(\theta) \end{array} \right),
\end{equation}
depend on $\theta$ strongly while the dependence on $v$ is less structured. 

We also note that keeping a non-zero diffusion matrix $\Sigma = 0.1I$ in (\ref{eq:PMP_forward}) leads to a desirable exploration of state space and, hence, to a stabilization/regularization of the control problem. However, a wider exploration of state space implies the need for larger ensemble sizes. Hence a balance of computational cost and state space exploration via added diffusion is required in practice and needs to be explored further.

%
\subsection{Lorenz-63: Infinite horizon control} \label{sec:num_L63}
%


\begin{figure}
\centerline{\includegraphics[width=0.9\textwidth]{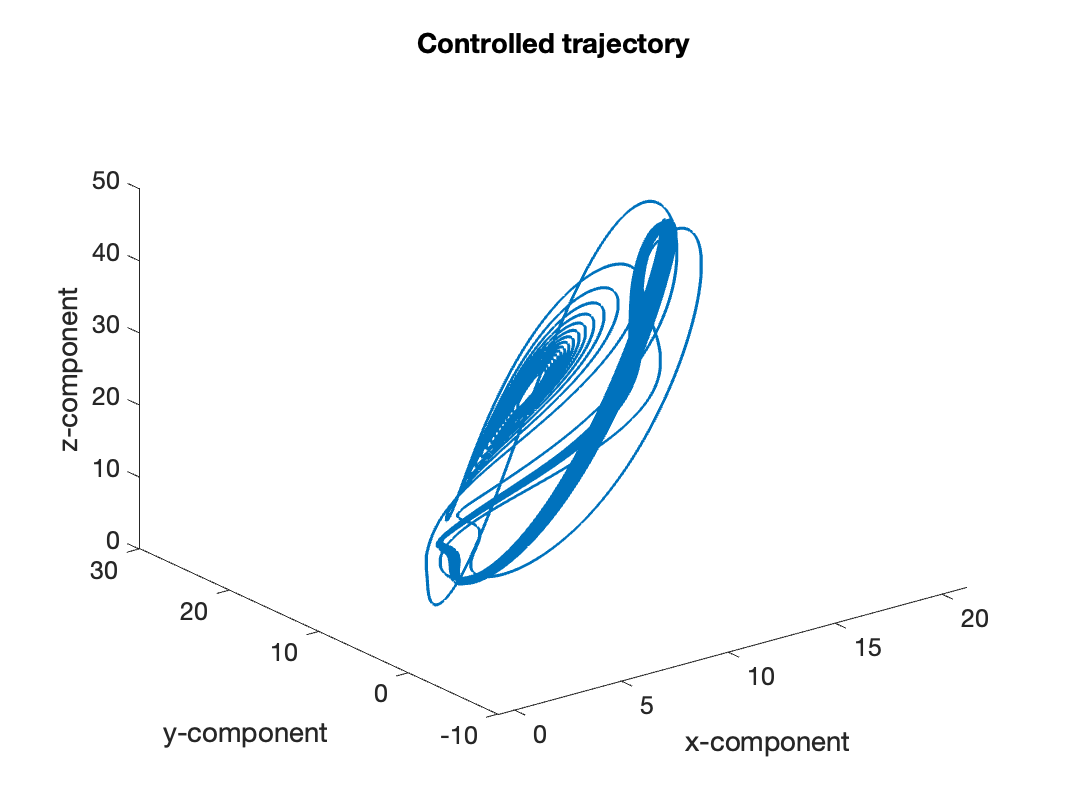}}
\caption{Controlled Lorenz-63 model: Displayed is the
three-dimensional trajectory of the  particle $\{X_t^{(i)}\}$ with $i=1$ over the time interval $t \in [0,100]$. After an initial transient, the trajectory enters an apparently quasi-periodic attractor.} 
\label{fig3}
\end{figure}

We follow \cite{KK24} and consider an optimal control problem for the Lorenz-63 model
\begin{equation} \label{eq:l63}
\dot{X}_t = b(X_t) + G U_t, \qquad b(x) = \left( \begin{array}{l} \sigma({\rm y}-{\rm x})\\
-{\rm x}{\rm z} + r{\rm x}-{\rm y}\\ {\rm x}{\rm y}-b{\rm z} \end{array}\right),
\end{equation}
in the state variable $x = ({\rm x},{\rm y},{\rm z})^{\rm T} \in \mathbb{R}^3$ and with
parameters $\sigma =10$, $r=28$, and $b = 8/3$ \cite{lorenz1963deterministic}. The scalar-valued control $U_t$ acts only on the ${\rm x}$-component; {\it i.e.}~$G = (1,0,0)^{\rm T}$. The control problem is to prevent ${\rm X}_t$ from taking negative values, and the authors of \cite{KK24} introduce the running cost 
\begin{equation} \label{eq:running cost L63}
c(x) = \frac{\alpha}{2} \left( \min({\rm x},0) \right)^2,
\end{equation}
which penalizes negative values of ${\rm X}_t$. In our experiments, we have followed this cost function and have set $\alpha = 5000$ and $R=1$
in (\ref{eq:norm}). Contrary to the setting of \cite{KK24}, we have implemented a discounted cost function (\ref{eq:costJ infinite}) with $\gamma = 1$ instead of a receding horizon NMPC approach \cite{Allgower,MPC}. 

We have implemented the Schr\"odinger bridge based McKean--Pontryagin formulation (\ref{eq:forward formulation}) with $\varepsilon = 0.1$, step-size $\Delta t = 0.001$, and $M=100$ particles. We again added artificial diffusion (viscosity) via $\Sigma = 0.5I$ in (\ref{eq:PMP_forward}). The parameter $\delta$ in (\ref{eq:NWKR}) is set to $\delta = 1$ and ${\rm b}_t = \dot{X}_t$ in (\ref{eq:NWKR_FD}). The initial momenta are $P_0^{(i)} = 0$.

\begin{figure}
\includegraphics[width=0.9\textwidth]{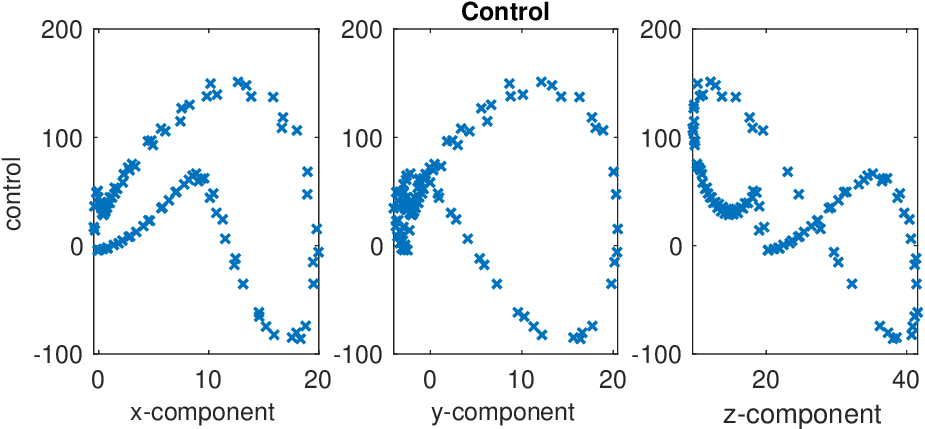}
\caption{Controlled Lorenz-63 model: Displayed is the dependence of the computed controls $U_T^{(i)}$ as a function of the three components of $X_T^{(i)}$ at final time $T = 100$ for $i=1,\ldots,M$, $M=100$.
} 
\label{fig2}
\end{figure}

The results can be found in Figures \ref{fig2} and \ref{fig3}. The impact of the control can be clearly detected. Due to the penalty induced by the
running cost (\ref{eq:running cost L63}), it is found that the controlled trajectories satisfy ${\rm X}_t \ge -0.5$. We also find that the control eliminates the chaotic nature of the Lorenz-63 model. The initial transition towards an apparently quasi-periodic attractor is due to the chosen initialization 
\begin{equation*}
X_0^{(i)} = (7.8590,7.1136,27.2293)^{\rm T} + \Xi_0^{(i)}, \qquad \Xi_0^{(i)} \sim {\rm N}(0,0.1I),
\end{equation*}
of the Lorenz-63 problem. Our results are qualitatively in agreement with the findings of \cite{KK24} which are based on a combination of data assimilation \cite{reich2015probabilistic,Evensenetal2022} and NMPC \cite{Allgower,MPC}. 

%
\subsection{Inverted pendulum: Receding horizon control} \label{sec:num_IP2}
%

We return to the inverted pendulum formulation (\ref{eq:pendulum}) from Section \ref{sec:num_IP1} but now in the receding horizon setting of NMPC. A control
is computed over receding windows of length $T=0.2$ by implementing (\ref{eq:linear variational}) with running cost given by (\ref{eq:running cost IP})
with $\alpha = 30$ and terminal cost
\begin{equation*}
    f(x) = \frac{\alpha_T}{2}(\theta-\pi)^2
\end{equation*}
with $\alpha_T = 200$. The computed closed loop control of the form
\begin{equation*}
    u_t(x) = -G(\theta)^{\rm T}(A_tx + c_t);
\end{equation*}
compare (\ref{eq:control law pendulum}) and (\ref{eq:linear regression}), is then applied to (\ref{eq:pendulum}) over a prediction interval of length $\Delta \tau = 0.02$ with diffusion (Brownian motion) added and $\Sigma =0.01I$. All equations are numerically integrated with step-size $\Delta t = 0.002$ and ensemble size $M=100$. The initial ensemble is drawn from the Gaussian ${\rm N}((2,0)^{\rm T}, 0.01I)$.

\begin{figure}
\centerline{\includegraphics[width=0.9\textwidth]{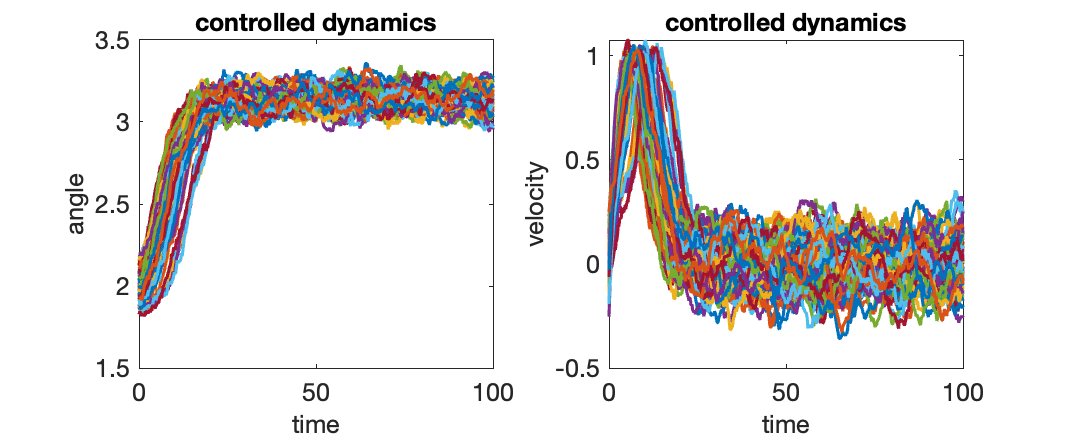}}
\caption{NMPC for inverted pendulum: Displayed are the time evolution of the ensemble of angles and velocities as a function of time. Starting from an ensemble centered at $(2,0)^{\rm T}$, the controlled ensemble eventually samples the vicinity of the unstable equilibrium $(\pi,0)^{\rm T}$. The magnitude of the fluctuations depends on the magnitude of the added diffusivity. Those fluctuations decrease as $\Sigma \to 0$.
} 
\label{fig4}
\end{figure}

We initialize each solution of (\ref{eq:linear variational}) with the positions from the last prediction step and the momenta from the previous optimization horizon and then solve  (\ref{eq:linear variational}) forward and backward once over a time window of length $T=0.2$. The numerical results can be found in Figure \ref{fig4}. We find that the receding horizon control formulation of the McKean--Pontryagin formulation is able to stabilize the unstable equilibrium $(\pi,0)^{\rm T}$. The remaining fluctuations in the ensemble are due to the non-vanishing diffusion which guaranties an exploration of state space in the vicinity of $(\pi,0)^{\rm T}$. Again, an optimal choice of $\Sigma$ will depend on the formulation of the particular problem. 

It should be noted that the variational formulation (\ref{eq:linear variational})
leads to local and linear approximations of the control, while the approach used in Section \ref{sec:num_IP1} delivers globally valid and nonlinear controls for the inverted pendulum. At the same time, the local approximation only requires the computation of covariance matrices instead of Schr\"odinger bridges and kernel regression, which suits the application to NMPC.

%
\subsection{Lorenz-96: Infinite horizon control} \label{sec:num_L96}
%

We test the McKean--Pontryagin minimum principle on the spatially extended $d_x=40$ dimensional Lorenz-96 model
\begin{align} \label{lorenz96}
\frac{{\rm d}x_l}{{\rm d}t} (t)=(x_{l+1}(t)-x_{l-2}(t))x_{l-1}(t)-x_l(t)  +10,\hspace{0.4cm}l=1,...,40
\end{align}
subject to periodic boundary conditions \cite{lorenz96}. The uncontrolled dynamics of the Lorenz-96 system leads to rapidly changing and highly chaotic solutions. The control task is to stabilize $x_{\rm ref} = 2$ with $G=I$; {\it i.e.}, the control $U_t\in \mathbb{R}^{40}$ acts on all components of the state variable $X_t = (x_1(t),\ldots,x_{40}(t))^{\rm T}$. We apply the infinite horizon formulation (\ref{eq:PMP_forward}) with running cost
\begin{equation}
c(x) = \frac{\rho}{2} \|x-x_{\rm ref}{\rm 1}_{40}\|^2,\quad 
\qquad {\rm 1}_{40} = (1,1,\ldots,1)^{\rm T} \in \mathbb{R}^{40},
\end{equation}
$\rho = 1000$, and discount factor $\gamma = 5$. The McKean--Pontryagin formulation (\ref{eq:PMP_forward}) is implemented with a $\Sigma = 0$ and linear approximation 
\begin{equation}
D_x^2 \phi_t(X_t)\dot{X}_t \approx {\rm C}_t^{px}({\rm C}_t^{xx})^{-1}\dot{X}_t.
\end{equation}
The ensemble size is set to $M=10$. Thus, $M\ll d_x$ and covariance localization is necessary. Covariance B-localization of ${\rm C}_t^{px}({\rm C}_t^{xx})^{-1}$ is based on the Gaspari--Cohn kernel \cite{gaspari99} with radius $r_{\rm loc} = 8$ \cite{reich2015probabilistic}. We also tested smaller localization radii, which led to comparable results. 

\begin{figure}
\centerline{\includegraphics[width=0.9\textwidth]{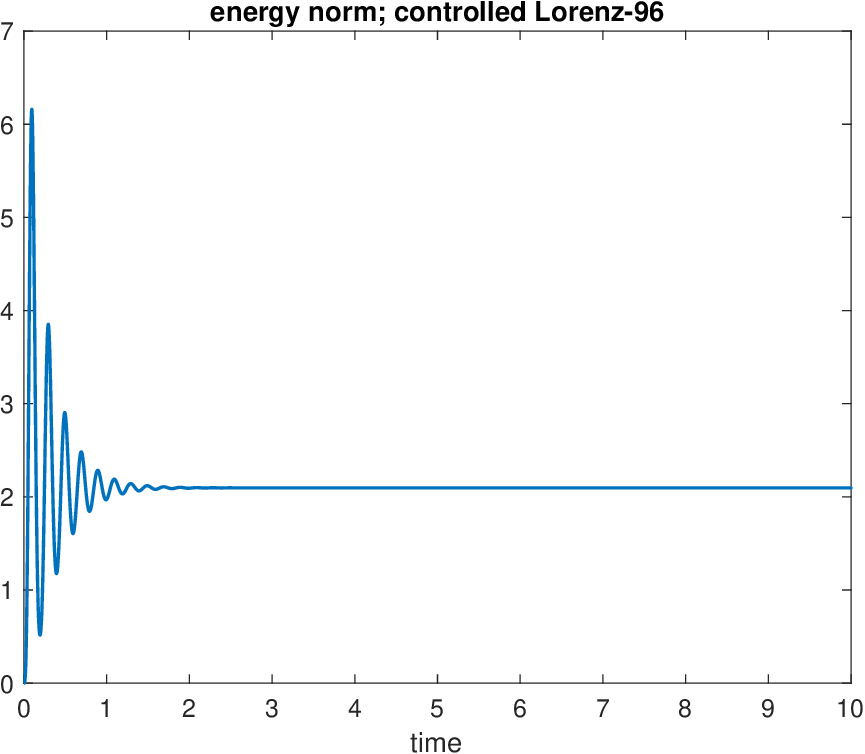}}
\caption{Top panel: Time evolution of the running cost (\ref{eq:energy norm}). It can be seen that the controlled dynamics stabilizes about the desired reference solution $x_{\rm ref}=2$. Bottom panel: The variance $\sigma(t_n)$ of the controlled ensemble $\{X_t^{(j)}\}$ converges to zero.
} 
\label{fig5}
\end{figure}

The system is initialized at $X_0 \sim {\rm N}(0,0.01I)$, $P_0 = 0$ and integrated up to time $T=5$. The step-size is $\Delta t = 0.001$. The time evolution of the running cost
\begin{equation} \label{eq:energy norm}
    c(t_n) = \frac{\rho}{2} \biggl\|\frac{1}{M} \sum_{i=1}^M X_{t_n}^{(i)}-x_{\rm ref}{\rm 1}_{40}\bigg\|^2
\end{equation}
along the ensemble mean can be found in Figure \ref{fig5} for $n = 0,\ldots,T/\Delta t$. The impact of control can be clearly seen and confirms that $x_{\rm ref} =2$ is effectively stabilized while uncontrolled solutions would develop into chaotic dynamics. We also show the ensemble variance $\sigma(t_n) = \mbox{trace} \,{\rm C}_{t_n}^{xx}$ in \ref{fig5}, which keeps decaying with increasing $t$.

%
\section{Conclusions}
%

A deterministic mean-field formulation of the PMP has been presented. Similar to the well-known forward-backward SDE formulation of the stochastic PMP \cite{Carmona} and its numerical implementation, a regression problem (\ref{eq:Hamiltonian_Lagrange}c) arises when discretized as an interacting particle system. Several numerical implementations have been discussed. Contrary to its stochastic counterpart, the McKean--Pontryagin minimum principle leads to smooth particle trajectories which should help with reducing the variance in the involved approximations. The freedom in the choice of the gauge variable $\beta_t$ allows for a decoupling of the forward and backward ODEs in $X_t$ and $P_t$, respectively. Furthermore, a forward simulation of the coupled Hamiltonian ODEs is possible in case of infinite horizon optimal control problems. Numerical results have been presented for a controlled one-dimensional diffusion process in a double well potential, a controlled inverted pendulum, a controlled Lorenz-63 system, and a controlled Lorenz-96 system. An application of the proposed methodology to NMPC \cite{MPC,Allgower} and a variational approximation of the McKean--Pontryagin formulation have also been numerically investigated. Although the latter experiments are qualitative in nature and demonstrate the broad scope of the proposed methodology, the first example also provides quantitative results in terms of achieved control values (\ref{eq:costJ}). 

An important area of future work is provided by extensions to partially observed optimal control problems \cite{ASTROM1965,Bensoussan92,PM13,9811560,R25a}. Here, a combination of the control formulation (\ref{eq:linear variational}) with the ensemble Kalman filter \cite{Evensenetal2022,CRS22} appears as a particularly promising direction, which leads to the appearance of mean-field variables in both the cost function (\ref{eq:costJ infinite}) as well as the forward SDE (\ref{eq:FSDE}). See \cite{R25c} for first results. Another promising application of our McKean--Pontryagin formulation is to generative modeling \cite{song2021scorebased}, which can be viewed as a stochastic optimal control problem \cite{BRU22}.

\medskip

\noindent
{\bf Acknowledgment.} This work has been funded by Deutsche Forschungsgemeinschaft (DFG) - Project-ID 318763901 - SFB1294.

\bibliographystyle{plainurl}
%
\bibliography{bib-database}
%

\end{document}